\documentclass[12pt]{article}
\usepackage{graphicx,amssymb,dsfont}
\usepackage[bookmarks=true]{hyperref}
\title{
Hydrodynamics and hydrostatics for a class of asymmetric particle
systems with open boundaries}
\author{
C. Bahadoran
\thanks{
Laboratoire de Math\'ematiques, Universit\'e Clermont-Ferrand 2,
F-63177 Aubi\`{e}re. e-mail:
Christophe.Bahadoran@math.univ-bpclermont.fr } }
\date{}
\newcommand{\dsp}{\displaystyle}
\newcommand{\bd}{\begin{displaymath}}
\newcommand{\be}{\begin{equation}}
\newcommand{\ba}{\begin{array}}
\newcommand{\ed}{\end{displaymath}}
\newcommand{\ee}{\end{equation}}
\newcommand{\ea}{\end{array}}
\newcommand{\espace}{\mbox{ }}

\newcommand{\Exp}{{\rm I\hspace{-0.8mm}E}}
\newcommand{\indicator}[1]{\mathds{1}_{#1}}

\newcommand{\abs}[1]{\left|#1\right|}
\def\N{\mathbb{N}}
\def\Z{\mathbb{Z}}
\def\R{\mathbb{R}}

\newcommand{\dt}{\partial_t}

\newcommand{\cspace}{\mathbb{E}}
\newcommand{\gen}{L}
\newcommand{\densities}{{\mathcal R}}

\newcommand{\eqref}[1]{(\ref{#1})}

\newtheorem{theorem}{Theorem}[section]

\newtheorem{proposition}{Proposition}[section]

\newtheorem{lemma}{Lemma}[section]
\newtheorem{corollary}{Corollary}[section]

\newenvironment{proof}[2]{\espace\\{\bf Proof of #1 \ref{#2}.}}{\hfill\mbox{$\square$}}

\begin{document}
\maketitle
\begin{abstract}
We consider attractive particle systems in $\Z^d$ with product
invariant measures.
We prove that when particles are restricted to a subset of $\Z^d$,
with birth and death dynamics at the boundaries, the hydrodynamic
limit is given by the unique entropy solution of a conservation law,
with boundary conditions in the sense of Bardos et al. (\cite{bln}).
For the hydrostatic limit between parallel hyperplanes, we prove a multidimensional version of the phase diagram conjectured in \cite{ps}, and
show that it is robust with respect to perturbations of the boundaries.
\end{abstract}
\section{Introduction}
\label{intro}
%
Stochastic lattice gases in contact with reservoirs are tractable and thus widely studied examples of nonequilibrium stationary states.
The derivation of the stationary macroscopic profile (hydrostatic limit) is a natural question in this context.
For diffusive systems, robust methods have been developed (see e.g.
\cite{els,els2,klo,lms}). The hydrostatic profile is the
stationary solution to the hydrodynamic equation (a possibly
nonlinear diffusion equation) with Dirichlet boundary conditions
imposed by reservoir densities. For instance, the symmetric simple exclusion
process exhibits a linear profile connecting these densities.
\\ \\
For driven lattice gases, the picture is different.
The asymmetric simple exclusion process with open boundaries was first introduced in \cite{lig2,lig3} as an intermediate tool for studying the process on $\Z^d$.
Its hydrostatic profile was determined by \cite{dehp} in the one-dimensional nearest-neighbor case. It consists of three phases with uniform bulk density:
low-density (LD) and high-density (HD) phases, where the bulk density is given by one of the boundaries,
and a maximum curent (MC) phase, where the bulk density is $1/2$.
Unlike in the diffusive setting, these profiles cannot satisfy both Dirichlet conditions if the reservoir densities are different.
LD and HD phases are separated by a coexistence line,
where the bulk state is a randomly located shock connecting the reservoir densities.
For more general models and currents, as well as higher space dimensions, mathematical results are missing.
In one space dimension, the number and nature of phases is expected to depend on the current-density function
through the following variational formula for the uniform bulk density (\cite{ps}):
\be\label{variational_ps}
{\rm argmin}_{[\lambda_a,\lambda_b]}f\mbox{ if }\lambda_a<\lambda_b,\quad
{\rm argmax}_{[\lambda_b,\lambda_a]}f\mbox{ if }\lambda_b<\lambda_a
\ee
where $\lambda_a$, $\lambda_b$ are the left and right reservoir densities, and
$f(\rho)$ is the current-density function. For the asymmetric exclusion process, $f(\rho)=\rho(1-\rho)$, and \eqref{variational_ps} yields the three  phases of \cite{dehp}.
For the KLS model (\cite{kls}), one obtains a seven-phase diagram, with two LD, two HD, two MC and a minimum current (mC) phase.
One outcome of this paper is to prove a multidimensional version of \eqref{variational_ps}, for a wide class of models including simple exclusion,
with arbitrarily many phases. The boundaries are basically parallel hyperplanes, but results are
somewhat robust with respect to perturbations of this geometry. The approach introduced here should be effective to treat more complex boundary-driven phase transitions, induced either by the domain geometry, or by two-species model like \cite{fri3}. These will be considered in future works. \\ \\
%
%
The key to our approach is to determine relevant boundary conditions in the scaling limit for asymmetric systems.
The hydrodynamic behavior of particle systems with open boundaries is
so far understood only in diffusive regimes (\cite{els2,lms,bemm,flm,mo}), where Dirichlet boundary conditions are relevant.
The celebrated result of \cite{rez} shows that, under Euler time scaling, the hydrodynamic limit of
attractive particle systems on $\Z^d$ with
product invariant measures, is given by the entropy solution to a scalar conservation law.
We extend this result to systems living in an open subset of $\R^d$.
%
%
We prove that the
hydrodynamic limit is the unique entropy solution to an initial-boundary
problem with BLN boundary conditions
introduced in \cite{bln}. Instead of fixing the boundary value like
Dirichlet conditions, these conditions impose a set of possible
boundary values depending on the boundary datum (\cite{dlf}).
Our proof uses a
generalized formulation (\cite{vov}) of the BLN boundary conditions
that does not explicitely involve a trace for the solution.
Producing this formulation from the microscopic boundary dynamics
involves adequate coupling of open systems.\\ \\
We next derive the hydrostatic profile and local equilibrium in a
domain lying between two parallel hyperplanes coupled to uniform
reservoirs. The result is a $d$-dimensional version of \cite{ps}: we show that the bulk density
is given by the variational formula of \cite{ps} applied to the normal projection of the flux. More generally, if the boundaries
are in some sense perturbations of hyperplanes, we show that, away from the perturbation, the bulk density
is the same as for hyperplanes,  regardless of the precise shape of the boundaries.
The hydrostatic limit follows from a uniqueness theorem that we establish for measure-valued
stationary entropy solutions with boundary conditions. Such a result implies (and is actually equivalent to)
asymptotic stability for entropy solutions with boundary conditions, a question studied so far
only for convex (\cite{lp,mt}) or bell-shaped (\cite{mv}) flux functions. We prove here such a result for general fluxes.\\ \\
We more generally expect BLN boundary conditions to arise in other models where convergence to the entropy solution is established on the whole space.
However, a proper microscopic treatment of the boundary remains to be found for models (e.g. \cite{sep,kls}) that are not attractive, or do not have explicit invariant measures. In the latter case, the ``natural'' definition of the boundary mechanism given here does not apply.
However, it is conjectured in \cite{gro} that the macroscopic behavior  of the boundary only depends on microscopic details of the boundary dynamics through an effective density.\\ \\
The paper is organized as follows. In section
\ref{section_framework} we define the framework and state the main results. In Section \ref{proof_hydro}, we establish the hydrodynamic (Theorems \ref{main}) and hydrostatic (Theorem \ref{main_2}) limits. The latter uses uniqueness of stationary solutions, proved in Section \ref{proof_main_3}, together with asymptotic stability
(Theorem \ref{main_3}). For simplicity, we prove the hydrodynamic and hydrostatic limits in detail for the Misanthrope's process,
formerly treated by \cite{rez} in $\Z^d$. However, as explained after Theorem \ref{main_2}, the scope of our results is best
illustrated by the exclusion process with overtaking, a model treated in Appendix \ref{appendix_kstep}. Despite computational
differences, our approach uses only attractiveness and product invariant measures, as suggested by the model-independent definition of the open dynamics
\eqref{def_genopen}
and associated coupling \eqref{def_couplegen}. The interested reader can find in an early version of this work (\cite{baba}) model-independent computations including both models in a unified framework.
\section{The framework and results}
\label{section_framework}
\textbf{Notation.} Let $\N=\{1,2,\ldots\}$ and
$\Z^+=\{0,1,\ldots\}$. For $x,y\in\R^d$ and $\delta>0$,
$\abs{x}=\max_{i=1,\ldots,d}\abs{x_i}$, $x\wedge y=\min(x,y)$,
$x\vee y=\max(x,y)$, $x^+=x\vee 0$, $x^-=-(x\wedge 0)$,
$B(x,\delta)=\{z\in\R^d:\,\abs{z-x}\leq\delta\}$. The $(d-1)$-dimensional Hausdorff measure on
$\R^d$ is denoted by ${\mathcal H}^{d-1}$.
The integral of a function $f$ against a measure $\mu$ will be denoted by $\int f\,d\mu$, $\mu(f)$ or $<\mu,f>$.\\ \\
We denote by $\cspace:=E^{\Z^d}$ the set of particle configurations
$\eta=(\eta(x):x\in\Z^d)$, where $E=[0,{\mathcal K}]\cap\Z^+$, and
$\mathcal K\in\N\cup\{+\infty\}$ is the maximum number of particles
per site. $\mathcal R:=[0,\mathcal K]\cap\R$ is the set of possible values
for the particle density.
A local function $f$ on $\cspace$ is a function that depends only on the restriction of $\eta$ to a finite subset of $\Z^d$.
Spatial shift on $\cspace$, defined by $\tau_x\eta(.)=\eta(x+.)$ for
$x\in\Z^d$, is extended to functions $f\in\R^\cspace$ by $\tau_x
f=f\circ\tau_x$, and to operators $L$
on $\R^\cspace$ by $\tau_x L=L\circ\tau_x$.\\ \\
The partial product order on $\cspace$, i.e. $\eta\leq\xi$ iff.
$\eta(x)\leq\xi(x)$ for all $x\in\Z^d$, induces (see e.g.
\cite{strassen, lig}) a partial stochastic order among probability
measures on $\cspace$: for two probability measures $\mu_1$ and
$\mu_2$ on $\cspace$, $\mu_1\leq\mu_2$ iff., equivalently: (a) For
every nondecreasing function $f$ on $\cspace$, $\int fd\mu_1\leq\int f
d\mu_2$; (b) There exists a probability measure $\tilde{\mu}$ on
$\cspace^2$, with marginals $\mu_1$ and $\mu_2$, supported on
$\{(\eta_1,\eta_2)\in{\cspace}^2:\,\eta_1\leq \eta_2\}$. \\ \\
%
%
%
%
%
%
%
\textbf{The model on $\Z^d$}. We consider here the Misanthrope's process introduced in \cite{coc}. Another model, the exclusion process with overtaking is described in Appendix \ref{appendix_kstep}.
Let $p(.)$ be a probability measure on $\Z^d$, with finite first moment,
satisfying the irreducibility assumption $\sum_{n\geq
1}[p^{*n}(.)+p^{*n}(-.)]>0$, where $*n$ denotes $n$-th power convolution.
%
%
Let $b(.,.)$ be a bounded function on
$E\times E$ such that $b(0,.)=0$, and either (if ${\mathcal
K}=+\infty$) $b(n,m)>0$ for $n>0$, or (if ${\mathcal K}<+\infty$)
$b(.,{\mathcal K})=0$ and $b(n,m)>0$ for $n>0$, $m<{\mathcal K}$. The
Misanthrope's process on $\Z^d$ is the Feller process on $\cspace$ with generator
\be\label{gen_misanthrope}
\gen f(\eta)=\sum_{{x,y}\in\Z}p(y-x)b(\eta(x),\eta(y))\left[f(\eta^{x,y})-f(\eta)\right]
\ee
for local functions $f:\cspace\to\R$, where $\eta^{x,y}$ is the new configuration resulting from $\eta$ after a particle
jumps from $x$ to $y$. Boundedness of $b(.,.)$, in the case $\mathcal K=+\infty$, is a technically simplifying but not necessary condition. It could be replaced by a Lipschitz
condition as in \cite{and}, in which case the process is constructed on a proper subset of $\cspace$.
We further assume that  $b(n,m)$ is nondecreasing with respect to $n$ and
nonincreasing with respect to $m$,
%
%
and satisfies the following algebraic conditions:
\be\label{condinvar} \left\{
\ba{lll}
\dsp\frac{b(n,m)}{b(m+1,n-1)} & = &
\dsp\frac{b(n,0)b(1,m)}{b(m+1,0)b(1,n-1)}\\ \\
b(n,m)-b(m,n) & = & b(n,0)-b(m,0)
\ea
\right.
\ee
Condition \eqref{condinvar} implies existence of a family  $(\nu_\rho,\,\rho\in\mathcal R)$ of product invariant measures whose
one-site marginal, denoted by $\theta_\rho$, is explicitely computable as follows. Let $q(n):=b(n,0)/b(1,n-1)$, that is a nondecreasing function, $q(n)!:=q(1)\cdots q(n)$. For $\beta\in[0,q(\infty))$ (the chemical potential), define the probability measure $\theta^\beta(n)=Z(\beta)^{-1}\beta^n/q(n)!$, where $Z(\beta)$ is a normalizing factor.
If $\mathcal K<+\infty$, it is extended by weak continuity to $\theta^\beta=\delta_\mathcal K$ for $\beta=q(\infty)$.
Its mean $R(\beta)$ is an increasing $C^\infty$ bijection from $[0,q(\infty)]\cap\R$ to
$\mathcal R$. We set $\theta_\rho:=\theta^{R^{-1}(\rho)}$.
The measure $\theta_\rho$, hence $\nu_\rho$, is stochastically nondecreasing with respect to $\rho$. By construction, $\nu_\rho[\eta(0)]=\rho$. Denoting by $\mathcal I$ the set of invariant measures for $\gen$,  and by $\mathcal S$ the set of shift-invariant probability distributions on $\cspace$, we have
\be\label{extremal_measures}
({\mathcal I}\cap{\mathcal S})_e=\{\nu_\rho:\,\rho\in\mathcal R\}
\ee
where index $e$ denotes extremal elements. The above model contains in particular the simple exclusion process
(where $\mathcal K=1$ and $b(n,m)=n(1-m)$) and zero-range process (where $\mathcal K=+\infty$ and $b(n,m)=g(n)$). Other examples can be found in \cite{ba1, bgrs}. \\ \\
%
%
%
%
%
\textbf{Open-boundary dynamics.}
Let $\Omega$ denote an open subset of $\R^d$ with locally finite perimeter, and $\Omega_N:=\{x\in\Z^d:\,x/N\in\Omega \}$ its lattice discretization, indexed by the scaling parameter $N\in\N$. Particles evolve in $\Omega_N$, while $\Z^d\backslash \Omega_N$ is a particle reservoir. Let $\lambda_N(.)$, the microscopic reservoir profile, be a $\mathcal R$-valued field on $\Z^d\backslash\Omega_N$. We assume that the sequence $(\lambda_N)$ is uniformly bounded, and is a lattice approximation, in a sense defined below (see \eqref{def_trace}), of a $\mathcal R$-valued  macroscopic reservoir profile, $\hat{\lambda}(.)\in L^\infty(\partial\Omega)$.\\ \\
We construct a Markov
generator $\gen_N$ on $\cspace_N:=E^{\Omega_N}$ as follows.
We decompose a particle configuration on $\Z^d$
as $\eta\oplus\overline{\eta}$, with $\eta\in
\cspace_N$ and $\overline{\eta}\in
\overline{\cspace}_N:=E^{\Z^d\backslash \Omega_N}$ (the reservoir state).
This decomposition depends on $N$, but for simplicity we omit this dependence in the notation.
We let $\overline{\eta}$ be distributed
according to a local equilibrium measure with density profile $\hat{\lambda}(.)$, that is a product measure
$\overline{\nu}_N$ on $\overline{\cspace}_N$, with one-site marginals given by
\be \label{outer_product}
\overline{\nu}_{N}(\{\overline{\eta}(x)=n\})=\theta_{\lambda_N(x)}(n),\quad
\forall x\in\Z^d\backslash \Omega_N,\,\forall n\in\Z^+ \ee
%
%
%
%
%
If $f(\eta)$ is a local function on $\cspace_N$, $Lf$ may depend on the inside plus reservoir configuration
$\eta\oplus\overline{\eta}$. We define a Markov generator by
\be \label{def_genopen} \gen_{N} f(\eta)=\int(\gen
f)(\eta\oplus\overline{\eta})d\overline{\nu}_{N}(\overline{\eta})\ee
%
%
%
Note that \eqref{def_genopen} is a model-independent way of defining  open-boundary dynamics given the dynamics on $\Z^d$,  so long as we have a family of product invariant measures for the latter.
%
%
For the Misanthrope's process we obtain the following explicit dynamics (see Subsection \ref{model_kstep} for the exclusion process with overtaking).
%
%
Define
\be \label{def_bbar} \overline{b}^+(\rho,n) =
\sum_{m}\theta_\rho(m)b(m,n),\quad %
\overline{b}^-(n,\rho) = \sum_{m}\theta_\rho(m)b(n,m) \ee
A jump in the bulk from $x\in\Omega_N$ to $y\in\Omega_N$ has the same rate $p(y-x)b(\eta(x),\eta(y))$ as in
\eqref{gen_misanthrope}. A reservoir particle at $x\in\Z^d\backslash\Omega_N$ jumps into the bulk to $y\in\Omega_N$, thus creating a particle at $y$, at rate $p(y-x)\overline{b}^+(\lambda_N(x),\eta(y))$; the total birth rate at $y$ is the sum of these contributions over $x$.
A particle in the bulk jumps from $x\in\Omega_N$ to the reservoir site $y\in\Z^d\backslash\Omega_N$, thus
removing a particle from $x$, at rate $p(y-x)\overline{b}^-(\eta(x),\lambda_N(y))$; the total death rate at $x$ is the sum of these contributions over $y$.\\ \\
{\em Example.}\label{page_example} Consider a one-dimensional totally asymmetric nearest-neighbor process: $p(1)=1$, $\Omega=(0,1)$, $\Omega_N=\{1,\ldots,N-1\}$, $\hat{\lambda}=\lambda_a\indicator{\{0\}}+\lambda_b\indicator{\{b\}}$,
$\lambda_N=\lambda_a\indicator{\{0\}}+\lambda_b\indicator{\{N\}}$.
For the simple exclusion process, $\overline{b}^+(\lambda,n)=\lambda(1-n)$, $\overline{b}^-(n,\lambda)=n(1-\lambda)$.
We then recover the entrance rate $\lambda_a[1-\eta(1)]$ and exit rate $\eta(N-1)[1-\lambda_b]$ of \cite{dehp}.\\ \\
{\em Remark.} For the zero-range process,  $b(n,m)$ depends only on $n$.
%
Thus, the entrance rate does not depend on $\eta$, i.e. the left reservoir behaves as a Poissonian source that does not see the bulk; and the exit rate does not depend on $\lambda_r$, i.e. the  bulk does not see the right reservoir.  This simplified boundary behavior has macroscopic counterparts, see remarks following theorems \ref{main} and \ref{main_2}.\\ \\
%
%
%
%
{\em Notational remark.} In the sequel, letters $\eta,\xi,\ldots$ will most often denote elements of $\cspace_N$. Occasionally, they will also denote elements
of $\cspace$. This will either be clear from the context, or mentioned explicitely.
%
%
%
%
%
%
%
\subsection{Hydrodynamic limit}
\label{hydro_open}
Consider the initial-boundary value problem on $\Omega$ for
\be \label{conservation_law} \dt\rho(t,x)+{\rm div}_x h(\rho(t,x))=0
\ee
where $h\in C^1(\mathcal R)$ and $||h'||_\infty<+\infty$.
The following definition of entropy solutions  is due to \cite{vov}, see also \cite{cf, ct} for related approaches. It differs
from the original definition of \cite{bln} and extensions thereof (\cite{ott, sz3, vas}) by the essential feature that it does not explicitely involve the trace
of the solution at the boundary, which makes it adapted to our problem. We say
$\rho(.,.)\in L^\infty((0,+\infty)\times\Omega)$ is an entropy solution to
\eqref{conservation_law} in $\Omega$, with initial datum $\rho_0\in
L^\infty(\Omega)$, and boundary datum $\hat{\lambda}\in L^\infty(\partial\Omega)$,
iff. there exists $M>0$ such that
\be \label{entropy_initial_vov} \ba{ll} &
\dsp\int_{(0,+\infty)\times\R^d}\left[
\partial_t\varphi(t,x)\phi(\rho(t,x))+\nabla_x\varphi(t,x).\psi(\rho(t,x))\right]dtdx
\\ + & \dsp
M\int_{(0,+\infty)\times\partial\Omega}\varphi(t,x)\phi(\hat{\lambda}(x))dt
\,d{\mathcal H}^{d-1}(x) +\int_{\Omega}\varphi(0,x)\phi(\rho_0(x))dx
\geq 0
 \ea\ee
for every  $0\leq \varphi\in C^1_K([0,+\infty)\times\R^d)$, and the family of Kru\v{z}kov entropy-flux pairs
(\cite{kru}) $(\phi,\psi)=(\phi_c^\pm,\psi_c^\pm)$ given for $c\in\mathcal R$ by
\be \label{def_semi} \ba{l} \phi^+_c(\rho)=(\rho-c)^+,\,
\psi^+_c(\rho)=1_{(0,+\infty)}(\rho-c)[h(\rho)-h(c)]\\
\phi^-_c(\rho)=(\rho-c)^-,\,
\psi^-_c(\rho)=-1_{(-\infty,0)}(\rho-c)[h(\rho)-h(c)] \ea \ee
%
%
More generally, let $\pi(dt,dx,d\rho)=dtdx\pi_{t,x}(d\rho)$ be a Young measure on $[0,+\infty)$, where $\pi$ is a weakly measurable mapping from $(0,+\infty)\times\Omega$
to $\mathcal P([0,+\infty))$ (the set of probability measures on $[0,+\infty))$. The set of Young measures is compact for the topology induced by the vague topology
for measures on $(0,+\infty)\times\R^d\times[0,+\infty)$. We say $\pi$ is bounded if there exists $C>0$ such that $\pi_{t,x}$ is supported a.e. on $[0,C]$.
In the spirit of \cite{dip},
the Young measure $\pi$ is called a mv solution of \eqref{conservation_law} with data $\rho_0\in
L^\infty(\Omega)$ and $\hat{\lambda}\in L^\infty(\partial\Omega)$, iff. the ``mv version'' of \eqref{entropy_initial_vov} holds, where the first integral is replaced by
\be\label{mv_version}
\int_{(0,+\infty)\times\R^d}\int_{[0,+\infty)}\left[
\partial_t\varphi(t,x)\phi(\rho)+\nabla_x\varphi(t,x).\psi(\rho)\right]\pi_{t,x}(d\rho)dxdt\ee
In particular, $\rho\in L^\infty((0,+\infty))\times\Omega$ is an entropy solution, iff. the Young measure $\pi_{t,x}=\delta_{\rho(t,x)}$ is a mv entropy solution.
It is clearly sufficient to require \eqref{entropy_initial_vov} or its mv version for a dense subset of values of $c$.
We can also define mv entropy solutions without an initial datum, by considering only $0\leq \varphi\in C^1_K((0,+\infty)\times\R^d)$, in which case the third integral in \eqref{entropy_initial_vov} is removed. In this case, the question of uniqueness is relevant only for stationary solutions (see Theorem \ref{main_3}).
For the latter, definition \eqref{entropy_initial_vov} reduces to considering only spatial test functions $\varphi\in C^1_K(\R^d)$ and removing time integrations.
For the Cauchy problem, we have the following result:
\begin{theorem}[\cite{vov}]\label{th_vovelle}
There exists a unique bounded mv entropy solution $\pi$ to \eqref{conservation_law} with data $\rho_0(.)$ and $\hat{\lambda}(.)$. This solution is of the form $\pi_{t,x}=\delta_{\rho(t,x)}$, where $\rho(t,x)$ is the unique entropy solution.
\end{theorem}
Let $\Delta=\sup\{|z|:\,z\in\Z^d,\,p(z)>0\}$.
We assume that the sequence $\lambda_N(.)$ has limiting trace $\hat{\lambda}(.)$ on $\partial\Omega$ in the following sense:
there exists a constant $C>0$ such that
\begin{eqnarray}
\limsup_{N\to\infty} & \dsp N^{1-d}\sum_{x\not\in\Omega_N,\,d(x,\Omega_N)\leq r}\varphi(x/N)f(\lambda_N(x)) & \nonumber\\
\leq & \dsp C(1+r)\int_{\partial\Omega}\varphi(x)f(\hat{\lambda}(x))\,d\mathcal H^{d-1}(x) & \label{def_trace}\end{eqnarray}
for every $r\in[0,\Delta]\cap\N$, $0\leq\varphi\in C^0_K(\R^d)$ and  $0\leq f\in C^0(\densities)$.
%
%
Condition \eqref{def_trace} holds for instance if
$\lambda_N(x)=\lambda(x/N)$, where $\lambda(.)$ is a continuous
function on $\R^d\backslash\Omega$ with trace $\hat{\lambda}(.)$ on
$\partial\Omega$.
Let $\mu^N$ be the probability distribution on $\cspace_N$ defined by
\be \label{loc_eq_2}
\mu^N(d\eta)=\bigotimes_{x\in\Omega_N}\nu_{\rho^N(x)}(d\eta(x))
\ee
where
$(\rho^N(x),\,N\in\N,x\in\Omega_N)$  is uniformly bounded,  $\mathcal R$-valued, and satisfies
\be \label{loc_eq_3} \lim_{N\rightarrow\infty}\int_I\abs{
\rho^N([Nx])-\rho(x) }dx=0 \ee
for every bounded Borel set $I\subset\Omega$. We assume that, for each $N\in\N$, $(\eta^N_t:\,t\geq 0)$ is a Markov process on $\cspace_N$
with generator $L_N$ and initial distribution $\mu^N$. Our first main result is the
\begin{theorem}
\label{main}
%
%
Let
$$\alpha^N_t(dx):=N^{-d}\sum_{y\in\Omega_N}\eta^N_{Nt}(y)\delta_{y/N}(dx)\in\mathcal M(\Omega)$$
where $\mathcal M(\Omega)$ denotes the set of Radon measures on $\Omega$ endowed with the topology of vague convergence. Then
for every $t>0$, $\alpha^N_t\to\alpha_t(dx):=\rho(t,x)dx$ in probability,
where $\rho(.,.)$ is the unique entropy solution to \eqref{conservation_law}
in $\Omega$ with initial datum $\rho_0(.)$ and boundary datum
$\hat{\lambda}(.)$, with flux function $h(.)$  given by
\be\label{current_macro}
h(\rho) = \int_\cspace j(\eta)\nu_\rho(d\eta)
\ee
where $j(\eta)$ is a microscopic flux function defined on $\cspace$ by
\begin{eqnarray} \label{current_misanthrope} j(\eta) & = & \sum_{z\in\Z^d}zp(z)b(\eta(0),\eta(z))
\end{eqnarray}
for Misanthrope's process, or \eqref{current_kstep} for the exclusion process with overtaking.
\end{theorem}
\textbf{Remark.} It follows from definition of $\nu_\rho$ that $h\in C^\infty(\mathcal R)$ and $||h'||_\infty<+\infty$.
Since $\nu_\rho$ is a shift-invariant product measure, for the Misanthrope's process,
%
%
$h(\rho)$ has constant direction $\gamma:=\sum_{z\in\Z^d}z p(z)$. Hence \eqref{conservation_law}
reduces to a family of uncoupled one-dimensional conservation laws.
The zero-range process is special because $h.\gamma$ is increasing. It follows (see \cite{dlf}) that an incoming BLN condition reduces to a Dirichlet condition, while
an outcoming BLN condition is void. This is a microscopic counterpart to the remark on page \pageref{page_example}.
For the model
of Appendix \ref{appendix_kstep}, the flux \eqref{flux_kstep} is genuinely multidimensional.
%
%
%
%
%
%
%
%
%
%
\subsection{Hydrostatic limit and stationary entropy solutions.}
%
%
Let $n\in\R^d$ be a unitary vector, $a,b\in\R$ with $a<b$, and
\be \label{hyperplanes} \Omega^n_{a,b}:=\{ x\in\R^d:\,n.x\in(a,b)
\} \ee
We consider \eqref{conservation_law} in $\Omega=\Omega^n_{a,b}$ with  boundary
datum
\be \label{boundary_datum_true_hyperplanes}
\hat{\lambda}(.)=\lambda_a \indicator{\{n.x=a\}}+\lambda_b \indicator{\{n.x=b\}}\ee
where $(\lambda_a,\lambda_b)\in{\mathcal R}^2$.
We want to study convergence to the stationary entropy solution, both for solutions of \eqref{conservation_law} and stationary distributions of $L_N$,
whenever this solution is unique. We first give a necessary condition for uniqueness.
Let  $f\in C^1(\mathcal R)$ be a real-valued  function such that $||f'||_\infty<+\infty$. Define
\begin{eqnarray}
\label{prime_1} \lambda^f_a & = &
\inf\{\lambda\leq\lambda_a:\,f\mbox{ is constant
on }[\lambda,\lambda_a]\}\\
\label{prime_2} \lambda^f_b & = &
\sup\{\lambda\geq\lambda_b:\,f\mbox{ is constant on
}[\lambda_b,\lambda]\}
\end{eqnarray}
if $\lambda_a\leq\lambda_b$, or
\begin{eqnarray}
\label{prime_3} \lambda^f_a & = &
\sup\{\lambda\geq\lambda_a:\,f\mbox{ is constant
on }[\lambda_a,\lambda]\}\\
\label{prime_4} \lambda^f_b & = &
\inf\{\lambda\leq\lambda_b:\,f\mbox{ is constant on
}[\lambda,\lambda_b]\}
\end{eqnarray}
if $\lambda_a\geq\lambda_b$.
For $\lambda_a\leq\lambda_b$ (resp. $\lambda_a\geq\lambda_b$), let
$\mathcal{M}_f(\lambda_a,\lambda_b)$ denote the set of minimizers
(resp. maximizers) of $f$ on
$[\lambda^f_a\wedge\lambda^f_b,\lambda^f_a\vee\lambda^f_b]$. If this set is reduced to a singleton, its unique element (which necessarily belongs to $[\lambda_a,\lambda_b]$) is denoted by $R_{f}(\lambda_a,\lambda_b)$.
For $\rho,\rho'\in\mathcal R$, we write $\rho=_f\rho'$ iff. $f$ is constant on the interval defined by $\rho,\rho'$, and $\rho\leq_f\rho'$ iff. $\rho\leq\rho'$ or $\rho=_f\rho'$.
\begin{proposition}\label{prop_stat}
Let $a=x_0<x_1<\ldots<x_n=b$, and $(\rho_0,\ldots,\rho_{n-1})$ be a nondecreasing (if $\lambda_a<\lambda_b$) or nonincreasing (if $\lambda_a>\lambda_b)$ sequence of elements
of $\mathcal{M}_{h(.).n}(\lambda_a,\lambda_b)$ in the sense of $\leq_{h(.).n}$. Then $\tilde{\rho}(x):={\rho}(n.x)$, where
\be\label{def_stat_sol}
{\rho}(x):=\sum_{k=0}^{n-1}\rho_k\indicator{(x_k,x_{k+1})}(x),
\ee
is a stationary entropy solution to \eqref{conservation_law} in $\Omega^n_{a,b}$ with boundary datum \eqref{boundary_datum_true_hyperplanes}.
\end{proposition}
\textbf{Remark.} For flux functions $h(.)$ obtained from  \eqref{current_macro}--\eqref{current_misanthrope} or \eqref{current_kstep}--\eqref{flux_kstep},
$h(.).n$ either has no flat segment, or is identically $0$. In the former case, $\lambda^f_\gamma=\lambda_\gamma$ for $\gamma\in\{a,b\}$, and $\leq_{h(.).n}$ reduces to $\leq$. In the latter case, any profile is a stationary solution.\\ \\
It follows from Proposition \ref{prop_stat} that
$\mathcal{M}_{h(.).n}(\lambda_a,\lambda_b)=\{R_{h(.).n}(\lambda_a,\lambda_b)\}$
is necessary for uniqueness of a stationary solution. We now assume
this condition satisfied, and establish uniqueness and convergence
to the unique solution. More generally, we consider a perturbation
$\Omega$ of $\Omega^n_{a,b}$
in the following sense: there exist $-\infty<a'\leq
a<b\leq b'<+\infty$ such that
\be
\label{perturbation}\Omega^n_{a,b}\subset\Omega\subset\Omega^n_{a',b'}
\ee
The boundary of $\Omega$ is a disjoint union of components
$\partial\Omega_\gamma$ for $\gamma\in\{a,b\}$, such that $n.x\leq
a$ (resp. $n.x\geq b$) on $\partial\Omega_a$ (resp. on
$\partial\Omega_b$).
%
%
Theorems \ref{main_3} and \ref{main_2} below state that, for the problem in $\Omega$ with boundary datum
\be \label{boundary_datum_hyperplanes}
\hat{\lambda}(.)=\lambda_a \indicator{\partial\Omega_a}+\lambda_b \indicator{\partial\Omega_b},\ee
uniqueness and convergence hold in $\Omega^n_{a,b}$. Thus, away from
the perturbation, the behavior of the system is the same as for
hyperplanes, regardless of the geometry of boundaries. Near
boundaries,  the density depends on the geometry. However, we prove
that it lies between the boundary datum and bulk density. This
implies that in LD or HD phases,
the bulk density extends up to the dominant boundary, regardless of the geometry of boundaries.\\ \\
In the following,  $u_n\to [\alpha,\beta]$ in $L^1_{\rm loc}(\overline{\Omega})$ means that, for every compact $K\subset\R^d$,
$\lim_{n\to\infty}\int_{K\cap\Omega}(u_n-\beta)^+dx=\lim_{n\to\infty}\int_{K\cap\Omega}(u_n-\alpha)^-dx=0$
\begin{theorem}
\label{main_3}
Assume $h\in C^1(\densities)$, $||h'||_\infty<+\infty$. Let $\rho(.,.)$, resp. $\pi$, be an entropy solution, resp. stationary mv entropy solution
%
%
%
%
to \eqref{conservation_law}, in $\Omega$ satisfying \eqref{perturbation},
with boundary datum \eqref{boundary_datum_hyperplanes}. Then the following hold, where $\rho^*=R_{h(.).n}(\lambda_a,\lambda_b)$,
and $I[\alpha,\beta]:=[\alpha\wedge\beta,\alpha\vee\beta]$:\\ \\
(i) $\pi_x(I[\lambda^f_a,\rho^*])=1$ on $\Omega\cap\Omega^n_{a',b}$, and
$\pi_x(I[\rho^*,\lambda^f_b])=1$ on $\Omega\cap\Omega^n_{a,b'}$.
In particular, $\pi_x=\delta_{\rho^*}$ a.e. on $\Omega^n_{a,b}$, which extends to $\Omega\cap\Omega^n_{a',b}$ (resp.
$\Omega\cap\Omega^n_{a,b'}$) if  $\rho^*=\lambda_a$ (resp. $\lambda_b$).\\ \\
(ii) As $t\to\infty$, $\rho(t,.)\to I[\lambda^f_a,\rho^*]$ in $L^1_{\rm loc}(\overline{\Omega\cap {\Omega}^n_{a',b}})$,
and $\rho(t,.)\to I[\rho^*,\lambda^f_b]$ in $L^1_{\rm loc}(\overline{\Omega\cap {\Omega}^n_{a,b'}})$.
In particular,
$\rho(t,.)\to \rho^*$ in $L^1_{\rm loc}(\overline{\Omega}^n_{a,b})$,
which extends to  $L^1_{\rm loc}(\overline{\Omega\cap {\Omega}^n_{a',b}})$, resp. $L^1_{\rm loc}(\overline{\Omega\cap {\Omega}^n_{a,b'}})$, if $\rho^*=\lambda_a$, resp. $\lambda_b$.
%
%
%
%
%
\end{theorem}
For the next result, we say a function $g:\cspace\to\R$ is Lipschitz, if there exists a finite $S\subset\Z^d$ and a constant $C>0$ such that
$|g(\eta)-g(\xi)|\leq C\sum_{x\in S}|\eta(x)-\xi(x)|$ for all $\eta,\xi\in\cspace$. This holds in particular for bounded local functions.
For $c\in\mathcal R$, we denote by $\nu^N_c$ the product measure on $\cspace_N$ under which $\eta(x)\sim\theta_c$ for every $x\in\Omega_N$.
\begin{theorem}
\label{main_2}
Let $h(.)$ be the flux function \eqref{current_macro}--\eqref{current_misanthrope} or \eqref{current_kstep}--\eqref{flux_kstep}.
Assume that, for each $N\in\N$, $\nu^N$ is an invariant measure for
$\gen_{N}$, and that there exists  $R\in\mathcal R$ such that $\nu^N\leq\nu^N_R$
for large enough $N\in\N$.
Let $\alpha^N(dx):=N^{-d}\sum_{y\in\Omega_N}\delta_{y/N}(dx)$.
Then, as $N\to\infty$: \\ \\
(i) the restriction of $\alpha^N$ to $\Omega^n_{a,b}$ converges  to $R_{h(.).n}(\lambda_a,\lambda_b)dx$
in probability under $\nu^N$ (hydrostatic limit); \\ \\
(ii) for every Lipschitz function $g:\cspace\to\R$,
$\overline{g}^N(x)  := \tau_{[Nx]}\nu^N[g(\eta)]
$
converges in $L^1_{\rm loc}(\overline{\Omega}^n_{a,b})$  to
$
\overline{g} := \nu_{R_{h(.).n}
(\lambda_a,\lambda_b)}[g(\eta)]
$
(local equilibrium).\\ \\
iii) If $R_{h(.).n}(\lambda_a,\lambda_b)=\lambda_a$, resp. $\lambda_b$, (i) and (ii) extend to
$\Omega\cap {\Omega}^n_{a',b}$, resp. $\Omega\cap {\Omega}^n_{a,b'}$.
\end{theorem}
\textbf{Remark.} The condition $\nu^N\leq\nu^N_R$ is a microscopic
counterpart of the boundedness assumption in the definition of stationary entropy solutions. If ${\mathcal K}<+\infty$, it is always satisfied by
$R={\mathcal K}$.
%
%
If ${\mathcal K}=+\infty$, there exists  a sequence of invariant
measures such that $\nu^N_r\leq\nu^N\leq\nu^N_R$, with $r=\min(\lambda_a,\lambda_b)$ and $R=\max(\lambda_a,\lambda_b)$.
Indeed, by Lemma \ref{monotone_coupling}, this holds for any subsequential limit of $M_t^N:=t^{-1}\int_0^t\nu^N_r
e^{s\gen_{N}}ds$ . In particular,
this is true if the invariant measure is
unique for $N$ large enough, which holds
in dimension one for any value of $\mathcal K$.\\ \\
%
%
%
For the Misanthrope's process, the flux $h$ has constant direction of the mean drift $\gamma$. Hence, the phase diagram depends on $n$ only through the sign of $n.\gamma$, and the bulk state is not determined if $n.\gamma=0$, where the flux is constantly parallel to the boundaries.
The zero-range process is special because $h.\gamma$ is increasing, thus the bulk density for $n.\gamma\neq 0$ is always that of the incoming boundary.
This is another  microscopic counterpart to the remark on page \pageref{page_example}.\\ \\
Note that the flux function \eqref{current_macro}--\eqref{current_misanthrope} is not explicit. The model of Appendix \ref{appendix_kstep} produces any $\R^d$-valued polynomial flux function (see \eqref{current_kstep}--\eqref{flux_kstep}). If $\mathcal K<+\infty$ and $h(.).n$ has $k$ local maxima and $k-1$ local minima, we obtain $k$ LD phases, $k$ HD phases, $k$ MC phases and $k-1$ mC phases.
%
%
%
%
%
\section{Proofs of Theorems \ref{main} and \ref{main_2}}
\label{proof_hydro}
%
%
%
We follow the scheme of \cite{rez} by deriving a mv  version
of entropy inequality \eqref{entropy_initial_vov} for the particle system, and using uniqueness results of Theorems \ref{th_vovelle} and \ref{main_3}.  The key novelty is to obtain the boundary term of \cite{vov}  from a proper coupling of  open systems.   Technically, the approach of \cite{rez} is slightly simplified by introducing a deterministic rather than random Young measure, and including the initial condition inside the entropy condition.
%
%
\subsection{Coupling open systems}
Let $c\in[0,\mathcal K]\cap\R$. We couple the $\eta$-process with reservoir profile $\lambda_N(.)$ with a $\xi$-process that has uniform reservoir profile $c$ (that is,
$(\xi_t)$ is a Feller process with generator similar to \eqref{def_genopen}, but $\lambda_N(.)$ replaced by the uniform profile $c$ on $\Z^d\backslash\Omega_N$.
%
%
The coupling is constructed in the same spirit as \eqref{def_genopen}. We start from the coupled $\tilde{\gen}$ generator for the dynamics on $\Z^d$
(see \cite{coc}),
\begin{eqnarray}\label{def_couplegen_bulk}
\tilde{L}f(\eta,\xi) & = &  \sum_{x,y\in\Z^d}b(\eta(x),\eta(y))\wedge b(\xi(x),\xi(y))\left[f(\eta^{x,y},\xi^{x,y})-f(\eta,\xi)\right]\nonumber\\
& + &  \sum_{x,y\in\Z^d}[b(\eta(x),\eta(y))-b(\xi(x),\xi(y))]^+\left[f(\eta^{x,y},\xi)-f(\eta,\xi)\right]\nonumber\\
& + &  \sum_{x,y\in\Z^d}[b(\eta(x),\eta(y))-b(\xi(x),\xi(y))]^-\left[f(\eta,\xi^{x,y})-f(\eta,\xi)\right]
\end{eqnarray}
where $(\eta,\xi)\in\cspace^2$.
Because of monotonicity assumptions on $b$,
this coupling preserves an initial order between $\eta$ and $\xi$.
Let $F_\rho(x)=\theta_\rho((-\infty,x])$  denote the c.d.f. of $\theta_\rho$, and $F_\rho^{-1}$ its generalized inverse. We define $\tilde{\theta}_{\rho,c}$ as the distribution
of $(F_\rho^{-1}(U),F_c^{-1}(U))$, where $U$ is a r.v. with uniform distribution on $[0,1]$. Since $\theta_\rho\leq \theta_c$ if $\rho\leq c$, $\tilde{\theta}_{\rho,c}$ is a monotone coupling of $\theta_\rho$ and $\theta_c$, i.e. it has these two measures as marginals, and if $\rho\leq c$, $\tilde{\theta}_{\rho,c}$ is supported on the set
of $(n,m)\in\Z^+$ such that $n\leq m$.
%
%
In the case $\mathcal K=1$, where the model reduces to simple exclusion, $\tilde{\theta}_{\rho,c}$ is the usual coupling of Bernoulli measures with parameters $\rho$ and $c$ , i.e.
\begin{eqnarray}\tilde{\theta}_{\rho,c}(0,0)  =  (1-\rho)\wedge(1-c) & , &
\tilde{\theta}_{\rho,c}(1,0)=(\rho-c)^+\nonumber\\
\tilde{\theta}_{\rho,c}(0,1) = (c-\rho)^+ & , &
\tilde{\theta}_{\rho,c}(0,0)=\rho\wedge c\label{bernoulli}
\end{eqnarray}
Let $\widetilde{\overline{\nu}}_N$ denote the product measure on $\overline{\cspace}_N^2$ with one-site
marginals given by
\be \label{outer_product_tilde}
\widetilde{\overline{\nu}}_{N}(\{\overline{\eta}(x)=n\})=\tilde{\theta}_{\lambda_N(x),c}(n),\quad
\forall x\in\Z^d\backslash \Omega_N,\,\forall n\in\Z^+ \ee
We define the coupled generator on $\cspace_N^2$ in a model-independent way by
\be\label{def_couplegen}
\tilde{\gen}_{N,c} f(\eta,\xi):=\int_{\overline{\cspace}_N^2}\tilde{L}f(\eta\oplus\overline{\eta},\xi\oplus\overline{\xi})d\widetilde{\overline{\nu}}_N(\overline{\eta},\overline{\xi})
\ee
where $(\eta,\xi)\in\cspace_N^2$ and $(\overline{\eta},\overline{\xi})\in\overline{\cspace}_N^2$.
More explicitely, for the Misanthrope's process, we obtain
\be\label{decomp_couple_open}
\tilde{\gen}_{N,c}=\tilde{\gen}_N^0+\tilde{\gen}_{N,c}^++\tilde{\gen}_{N,c}^-
\ee
where $\tilde{\gen}^N_0$ is similar to $\tilde{\gen}$ in \eqref{def_couplegen_bulk}, but with $x,y$ restricted to $\Omega_N$, while
\begin{eqnarray}
\tilde{\gen}_{N,c}^+ f(\eta,\xi)  =   \sum_{x\not\in\Omega_N,\,y\in\Omega_N}p(y-x)\overline{b}^+(\lambda_N(x),\eta(y))\wedge \overline{b}^+(c,\xi(y))\left[f(\eta+\delta_y,\xi+\delta_y)-f(\eta,\xi)\right] & & \nonumber\\
+  \sum_{x\not\in\Omega_N,\,y\in\Omega_N}p(y-x)[\overline{b}^+(\lambda_N(x),\eta(y))- \overline{b}^+(c,\xi(y))]^+\left[f(\eta+\delta_y,\xi)-f(\eta,\xi)\right] & & \nonumber\\
+ \sum_{x\not\in\Omega_N,\,y\in\Omega_N}p(y-x)[\overline{b}^+(\lambda_N(x),\eta(y))-\overline{b}^+(c,\xi(y))]^-\left[f(\eta,\xi+\delta_y)-f(\eta,\xi)\right] & & \label{def_couplegen_plus}
\end{eqnarray}
\begin{eqnarray}
\tilde{\gen}_{N,c}^- f(\eta,\xi) =  \sum_{x\in\Omega_N,\,y\not\in\Omega_N}p(y-x)\overline{b}^-(\eta(x),\lambda_N(y))\wedge \overline{b}^-(\xi(x),c)\left[f(\eta-\delta_x,\xi-\delta_x)-f(\eta,\xi)\right] & & \nonumber\\
+  \sum_{x\in\Omega_N,\,y\not\in\Omega_N}p(y-x)[\overline{b}^-(\eta(x),\lambda_N(y))- \overline{b}^-(\xi(x),c)]^+\left[f(\eta-\delta_x,\xi)-f(\eta,\xi)\right] & & \nonumber\\
+  \sum_{x\in\Omega_N,\,y\not\in\Omega_N}p(y-x)[\overline{b}^-(\eta(x),\lambda_N(y))-\overline{b}^-(\xi(x),c)]^- \left[f(\eta,\xi-\delta_x)-f(\eta,\xi)\right] & & \label{def_couplegen_min}
\end{eqnarray}
%
%
%
%
%
In the sequel, $(\eta^N_t,\xi^N_t)$ denotes a coupled process with
generator $\tilde{L}_{N,c}$.
Expectation  with respect to this process will be denoted by $\Exp_{N,c}$ when we want to emphasize dependence on $c$, otherwise simply by $\Exp$.
We assume that $\xi^N_0\sim\nu^N_c$ and that, for some $R\in\mathcal R$, the distribution of $\eta^N_0$ is dominated stochastically by
$\nu^N_R$ for all $N\in\N$. The latter is automaticallly fulfilled with $R=\mathcal K$ if $\mathcal K<+\infty$.
\begin{lemma}
\label{equilibrium_property}
For every $t\geq 0$, $\xi^N_t$ has distribution $\nu^N_c$.
\end{lemma}
\begin{proof}{lemma}{equilibrium_property}
Let $\gen_{N,c}$ denote the generator obtained from \eqref{def_genopen}
when the reservoir profile $\lambda_N(.)$ has uniform value $c$.
Let $f(\eta)$ be a local function on $\cspace_N$.
Since $\nu_c$ is invariant for $\gen$,
$$
\int_{\cspace_N}\gen_{N,c}f(\eta)d\nu^N_{
c}(\eta)=\int_{\cspace}\gen
f(\eta\oplus\overline{\eta})d\nu_c(\eta\oplus
\overline{\eta})=0
$$
Hence $\nu^N_c$ is an invariant measure for \eqref{def_genopen}.
%
\end{proof}
\begin{lemma}
\label{monotone_coupling}
For every $t\geq 0$, $\eta^N_t$ is dominated stochastically by $\nu^N_\Lambda$, where
$\Lambda:=\max[R,\sup_{x\in\Z^d\backslash\Omega_N}\lambda_N(x)]$.
\end{lemma}
\begin{proof}{lemma}{monotone_coupling}
Consider a coupled process $(\eta^N_t,\zeta^N_t)$ with generator $\tilde{L}_{N,\Lambda}$, whose initial distribution is such that $\eta^N_0\leq\zeta^N_0$ and $\zeta^N_0\sim\nu^N_\Lambda$. The initial order is preserved by jumps within $\Omega_N$, since the coupling inside $\Omega_N$ is the same as in $\Z^d$. It is not hard to see that  the order is also preserved by the coupling of births or deaths, because $\overline{b}^+(\rho,n)$ and $\overline{b}^-(\rho,n)$, defined in \eqref{def_bbar} are respectively nondecreasing and nonincreasing in $\rho$.
\end{proof}
%
%
\subsection{Average entropy inequality}
To prove Theorem \ref{main}, we define a Young measure that we show satisfies \eqref{entropy_initial_vov}.
%
%
Let $\widetilde{\mathcal I}$ denote the set of invariant measures for \eqref{def_couplegen_bulk}, and
$\widetilde{\mathcal S}$ the set of shift-invariant probability distributions on $\cspace^2$. The arguments used for the proof of \eqref{extremal_measures} also establish the following result:
\begin{proposition}\label{invariant_coupled}(\cite{and, coc, gui, lig})
$$(\widetilde{\mathcal I}\cap\widetilde{\mathcal S})_e=\{\widetilde{\nu}_{\rho,\rho'}:\,(\rho,\rho')\in{\mathcal R}^2\}$$
where $\widetilde{\nu}_{\rho,\rho'}(d\eta,d\xi)$ is a probability distribution on $\cspace^2$ with marginals $\nu_\rho(d\eta)$ and $\nu_\rho'(d\xi)$,
and is supported on the set $\{\eta\leq\xi\}$ (resp. $\geq$)
if $\rho\leq\rho'$ (resp. $\geq$).
\end{proposition}
Let $\mathcal F$ denote the set of functions $F(t,x,\eta,\xi):[0,+\infty)\times\R^d\times\cspace^2\to\R$, uniformly continuous in $(t,x)$, with
(i) compact support: $F(t,x,\eta,\xi)=0$ if $(t,x)\not\in K$; (ii) Lipschitz property in $(\eta,\xi)$: $|F(t,x,\eta,\xi)-F(t,x,\eta',\xi')|\leq C\sum_{|z|\leq r}(|\eta(z)-\eta'(z)|+|\xi(z)-\xi'(z)|)$,
where $r\in\N$, $C>0$ are constants, and $K$ a compact subset of $\subset[0,+\infty)\times\R^d$.
\begin{proposition}\label{loc_eq}
Let $\mathcal S$ be a dense subset of $\mathcal R$.
For any  sequence of integers increasing to $+\infty$, there exists a subsequence and a Young measure $(t,x)\mapsto\pi_{t,x}$ supported a.e. on $[0,\Lambda]$,
such that, for every $F\in\mathcal F$
and $c\in\mathcal S$,
\begin{eqnarray}
\lim_{N\to\infty} & \Exp
_{N,c}
\left\{
\dsp\int_0^{+\infty}N^{-d}\sum_{x\in\Omega_N:\,d(x,\Z^d\setminus\Omega_N)>r}F\left(t,\frac{x}{N},\tau_x\eta^N_{Nt},\tau_x\xi^N_{Nt}\right)dt
\right\} & \nonumber\\
= & \dsp\int_{0}^{+\infty}\int_{\Omega}\int_{[0,+\infty)}\left<F(t,x,.,.)\right>_{\rho,c}\pi_{t,x}(d\rho)dxdt & \label{convergence_young}
\end{eqnarray}
holds along this subsequence,
where $<>_{\rho,c}$ denotes expectation w.r.t. $\widetilde{\nu}_{\rho,c}$.
\end{proposition}
Let $0\leq \varphi\in C^1_K([0,+\infty)\times\R^d)$.
Define
\be\label{def_fpm}
\Phi^\pm_N(t,\eta,\xi):=N^{-d}\sum_{x\in\Omega_N}\varphi(t,x/N)(\eta(x)-\xi(x))^{\pm}
\ee
Then,
\be
\Exp\left\{
\Phi^\pm_N(0,\eta^N_{0},\xi^N_{0})+ \int_0^{+\infty}
[\partial_t+
N\tilde{L}_{N,c}]\Phi^\pm_{N}(t,\eta^N_{Nt},\xi^N_{Nt})dt\right\}=0
\label{def_martingale}
\ee
%
%
%
%
%
We denote respectively by  $\mathcal K^\pm_{\partial\Omega}(\varphi)$ and $\mathcal K^\pm_0(\varphi)$, the second and third integrals in \eqref{entropy_initial_vov}, and by $\mathcal K^\pm_\Omega(\varphi,\pi)$ the integral in \eqref{mv_version}. The main step in the proof of Theorem \ref{main} is the following:
\begin{proposition}\label{micro_vovelle}
Let $\pi$ be a Young measure given by Proposition \ref{loc_eq}.
Then \eqref{replacement_boundary}--\eqref{replacement_bulk} below hold, where
the constant $M>0$ depends only on $p(.)$ and $b(.,.)$. In addition, \eqref{replacement_0} holds if the initial distribution of $(\eta^N_0,\xi^N_0)$ is the product measure $\tilde{\mu}^{N,c}$ on $\cspace_N^2$ whose marginal at site $x\in\Omega_N$ is
$\tilde{\theta}_{\rho^N(x),c}$.
\begin{eqnarray}
\lim_{N\to\infty} \Exp\Phi^\pm(0,\eta^N_0,\xi^N_0) & = & \mathcal K^\pm_0(\varphi)\label{replacement_0}\\
\limsup_{N\to\infty} \Exp\left\{
\int_0^{+\infty}
N[\tilde{L}^+_{N,c}+\tilde{L}^-_{N,c}]\Phi^\pm_{N}(t,\eta^N_{Nt},\xi^N_{Nt})dt
\right\} & \leq & M\mathcal K^\pm_{\partial\Omega}(\varphi)\label{replacement_boundary}\\
\limsup_{N\to\infty} \Exp\left\{
\int_0^{+\infty}
[\partial_t+
N\tilde{L}^0_{N,c}]\Phi^\pm_{N}(t,\eta^N_{Nt},\xi^N_{Nt})dt
\right\} & \leq & \mathcal K^\pm_{\Omega}(\varphi,\pi)\label{replacement_bulk}
\end{eqnarray}
\end{proposition}
From this and \eqref{def_martingale}, we immediately deduce the following:
\begin{corollary}\label{young_mv}
Any limiting Young measure $\pi$ in Proposition \ref{loc_eq} is a mv entropy solution to \eqref{conservation_law} on $\Omega$ with boundary datum $\hat{\lambda}(.)$.
If $(\eta^N_0,\xi^N_0)\sim\tilde{\mu}^{N,c}$ defined in Proposition \ref{micro_vovelle}, then $\pi$ has initial datum $\rho_0(.)$.
\end{corollary}
\textbf{Proof of Theorems \ref{main} and \ref{main_2}.}  Under
assumptions of Theorem \ref{main}, Corollary \ref{young_mv} and
Theorem \ref{th_vovelle} imply that the whole sequence in
\eqref{convergence_young} converges to the limit given by
$\pi_{t,x}=\delta_{\rho(t,x)}$, where $\rho(t,x)$ is the unique
entropy solution. To deduce convergence of the empirical measure is
a standard technical step (see e.g. \cite{rez,kla}). Similarly,
Corollary \ref{young_mv} combined with Theorem \ref{main_3} implies
Theorem \ref{main_2}.
\\ \\
\begin{proof}{proposition}{micro_vovelle}
\mbox{}\\ \\
{\em Proof of \eqref{replacement_0}.}
This follows from \eqref{loc_eq_3} since, by definition of $\tilde{\mu}^{N,c}$,
$$
\Exp\Phi_N^\pm(0,\eta^N_{0},\xi^N_{0})=N^{-d}\sum_{x\in\Z^d}\varphi\left(0,\frac{x}{N}\right)(\rho^N(x)-c)^\pm$$
{\em Proof of \eqref{replacement_boundary}.}
let $\tilde{\mathcal L}\in\{\tilde{\gen}^+_{N,c},\tilde{\gen}^-_{N,c}\}$. Then, for some constant $C>0$, depending only on $p(.)$ and $b(.,.)$,
\be\label{kruzkov_boundary}
N\tilde{\mathcal L}\Phi_N^{{\pm}}(t,\eta,\xi) \leq
C N^{1-d}\sum_{x\in\Omega_N}\varphi(t,x/N)\sum_{y\not\in\Omega_N}[p(y-x)+p(x-y)](\lambda_N(y)-c)^{{\pm}}
\ee
Indeed, consider for instance $\Phi^+_N$. The only terms in \eqref{def_couplegen_plus}--\eqref{def_couplegen_min} that produce a positive variation of $\Phi_N^+$
are those of \eqref{def_couplegen_plus} for which $\eta(y)\geq\xi(y)$ and a particle is created at $y$ in $\eta$ only, and those of \eqref{def_couplegen_min} for which $\eta(x)\geq\xi(x)$ and a particle is removed at $x$ from $\xi$ only.
Then, \eqref{kruzkov_boundary} is a consequence of the following
inequalities for $\eta(x)\geq\xi(x)$, $\eta(x)\geq\xi(y)$:
\begin{eqnarray*}
\{\overline{b}^+(\lambda_N(x),\eta(y))-\overline{b}^+(c,\xi(y))\}^+ & \leq & ||b||_\infty(\lambda_N(x)-c)^+\\
\{\overline{b}^-(\eta(x),\lambda_N(y))-\overline{b}^-(\xi(x),c)\}^- & \leq & ||b||_\infty(\lambda_N(y)-c)^+
\end{eqnarray*}
The above inequalities are consequences of monotonicity assumptions on $b(.,.)$ and
the fact that $\tilde{\theta}_{\lambda,c}$ is a monotone coupling of $\theta_\lambda$ and $\theta_c$.
By \eqref{def_trace} and the first moment assumption for $p(.)$, in the limit $N\to\infty$, the r.h.s. of \eqref{kruzkov_boundary} is bounded by a uniform constant times $\mathcal K^\pm_{\partial\Omega}(\varphi)$. \\ \\
{\em Proof of \eqref{replacement_bulk}.}
The same computation as in \cite[equation (3.2)]{rez} (with minor differences due to the boundary and finite-range cutoff) shows that, for every $r\in\N$,
%
%
%
%
\be\label{computation_rez}
N\tilde{\gen}^0_{N}\Phi^\pm_N(t,\eta,\xi)  \leq  N^{-d}\sum_{x\in\Omega_N:\,d(x,\Z^d\setminus\Omega_N)>r}\nabla_x\varphi(t,x/N).\tau_x\tilde{j}_r^\pm(\eta,\xi)+\delta_r+\varepsilon_N
%
%
%
%
\ee
where $\lim_{N\to\infty}\varepsilon_N=\lim_{r\to\infty}\delta_r=0$, and
\begin{eqnarray}
\tilde{j}_r^{+}(\eta,\xi) & = & \sum_{z\in\Z^d:\,|z|\leq r}zp(z)F_{0,z}(\eta,\xi)[b(\eta(0),\eta(z))-b(\xi(0),\xi(z))]\label{current_plus}\\
%
%
\tilde{j}_r^-(\eta,\xi) & = & \sum_{z\in\Z^d:\,|z|\leq r}zp(z)F_{0,z}(\xi,\eta)[b(\xi(0),\xi(z))-b(\eta(0),\eta(z))]\label{current_min}\\
F_{x,y}(\eta,\xi) &  =  & \indicator{\{\eta(x)\geq\xi(x),\,\eta(y)\geq\xi(y)\}}
%
%
%
\end{eqnarray}
In particular, $\tilde{j}^+_r(\eta,\xi)=j_r(\eta)-j_r(\xi)$ on $\{\eta\geq\xi\}$, resp. $0$ on $\{\eta\leq\xi\}$, where $j_r$ is defined as \eqref{current_misanthrope}
with the truncation $|z|\leq r$. Thus by Proposition \ref{invariant_coupled},
$$<(\eta(0)-\xi(0)^+>_{\rho,c}=(\rho-c)^+=\Phi_c^+(\rho),\quad
<\tilde{j}^+_r>_{\rho,c}=\indicator{\rho\geq c}[h_r(\rho)-h_r(c)]$$
where $h_r$ is defined from $j_r$ as in \eqref{current_macro}. Since $h_r\to h$ uniformly,
the result follows from Proposition \ref{loc_eq}.
\end{proof}
\mbox{}\\
\begin{proof}{proposition}{loc_eq}
It is enough to  prove existence of a Young measure $\pi$ satisfying \eqref{convergence_young} for a given $c\in\mathcal S$. Indeed, by diagonal extraction, we can then find
a common subsequence of $N\to\infty$ along which, for each $c\in\mathcal S$, there is a Young measure $\pi^c$, a priori depending on $c$. However, taking in \eqref{convergence_young} arbitrary
test functions not depending on $\xi$, shows  that $\pi^c$ does not depend on $c$. \\ \\
Given $a=(a_1,\ldots,a_d)$ and $b=(b_1,\ldots,b_d)$ in $\R^d$ such that $a\leq b$ componentwise, we set $[a,b):=\prod_{i=1}^d[a_i,b_i)$.
For $n\in\N$, $k=(k_1,\ldots,k_d)\in\Z^d$ and $l\in\N$, we set $x_{n,k}=k\,2^{-n}\in\R^d$, $B_{n,k}=[k\,2^{-n},(k+1)\,2^{-n})$, $t_{n,l}=l\,2^{-n}$ and $T_{n,l}=[l\,2^{-n},(l+1)\,2^{-n})$.
Let $\widetilde{\mu}^N_{t}$ denote the distribution of $(\eta^N_{Nt},\xi^N_{Nt})$, and
\be\label{average_measure}M^N_{n,k,l}:=2^{n}\int_{T_{n,l}}\frac{1}{N^d 2^{nd}}\sum_{x\in B_{n,k}}\tau_x\widetilde{\mu}^N_{t}dt\ee
Below we show that there exists a constant
$V>0$ with the following property: for each $(k,l)\in\Z^d\times\Z$ such that
\be\label{condition_box}B_{n,k}\subset\Omega,\quad d(B_{n,k},\R^d\setminus\Omega)> V.2^{-n}\ee
there is a subsequence along which
\be\label{averaged_measure} M^N_{n,k,l}\to\widetilde{\nu}_{n,k,l}:=\int_{\mathcal R}\widetilde{\nu}_{\rho,c}d\pi_{n,k,l}(\rho)\in\widetilde{\mathcal I}\cap\widetilde{\mathcal S}\ee
where  $\pi_{n,k,l}$ is a probability measure on $[0,\Lambda]$, and  $\to$ means convergence for Lipschitz functions on $\cspace^2$. By diagonal extraction, there is a common subsequence along which this convergence holds simultaneously
for all triples $(n,k,l)$.  Define the Young measure
$$ \pi^n_{t,x}:=\sum_{k\in\Z^d,\,l\in\Z}\pi_{n,k,l}\indicator{T_{n,l}}(t)\indicator{B_{n,k}}(x)$$
%
%
By compactness, there is a limiting (in vague topology) Young measure $\pi_{t,x}$ for $\pi^n_{t,x}$ along a subsequence of $n\to\infty$.
Since $\pi^n_{t,x}$ is supported a.e. on a fixed interval $[0,\Lambda]$, the same holds for $\pi_{t,x}$,
and vague convergence extends to $<F(t,x,.,.>_{\rho,c}$, even though it is not bounded.
The result follows by continuity assumption on $F$.\\ \\
{\em Proof of \eqref{average_measure}.} By Lemma \ref{monotone_coupling}, on each finite subset of $\Z^d$, $M^N_{n,k,l}$ has both marginals dominated by $\nu_{\Lambda\vee c}$.
Thus for fixed $n,k,l$ it is a tight sequence as $N\to\infty$, and for every subsequential limiting distribution, convergence holds on functions with at most linear growth.
%
%
We couple our process $(\eta^N_{t},\xi^N_{t})$ on the time interval $NT_{n,l}$ to a process $(\hat{\eta}^N_t,\hat{\xi}^N_t)$ on $\cspace^2$ with generator \eqref{def_couplegen_bulk}, so that in $(\eta^N_t,\hat{\eta}^N_t)$ and $(\xi^N_t,\hat{\xi}^N_t)$, jumps within
$\Omega_N$ are coupled as in \eqref{def_couplegen_bulk}.
%
%
At time $Nt_{n,l}$, the initial distribution of the coupling is
chosen so that $\hat{\eta}^N$ coincides with $\eta^N$ on $\Omega_N$
and has no particle outside, and similarly for $\hat{\xi}^N$. Denote
by $\hat{M}^N_{n,k,l}$ the measure defined by replacing
$\widetilde{\mu}^N_t$ in \eqref{averaged_measure} with the
distribution of $(\hat{\eta}^N_{Nt},\hat{\xi}^N_{Nt})$. Statement
\eqref{average_measure} for $\hat{M}^N_{n,k,l}$ follows as in the
proof of \cite[Theorem 3.1]{rez} from shift-invariance of
\eqref{def_couplegen_bulk}, and is then deduced for $M^N_{n,k,l}$
from the following lemma.
\end{proof}
\begin{lemma}
\label{mixed_comparison}
Under condition \eqref{condition_box},
$$\lim_{N\to\infty}\sup_{t\in T_{n,l}}\Exp N^{-d}\sum_{x\in N
B_{n,k}}\abs{\eta^N_{Nt}(x)-\hat{\eta}^N_{Nt}(x)}=0$$
and similarly for $(\xi^N_{Nt},\hat{\xi}^N_{Nt})$.
\end{lemma}
%
%
%
%
\begin{proof}{lemma}{mixed_comparison}
Let $\varepsilon>0$, $H_\varepsilon(.)=H_1(./\varepsilon)$, where $H_1$ is a smooth nondecreasing function on $\R$ such that $H_1(x)=0$ for $x<0$, $H_1(x)=1$ for $x\geq 1$.
For $t\in T_{n,l}$, let $\Phi_N=\Phi^+_N+\Phi^-_N$, see \eqref{def_fpm}, with
$\varphi(t,x)=H_\varepsilon(R-V(t-t_{n,l})-|x-x_{n,k}|_\varepsilon)$, where $V:=\sum_{z\in\Z^d}|z|p(z)$ and $|x|_\varepsilon:=(|x|^2+\varepsilon^2)^{1/2}$, with $R>0$ and $\varepsilon>0$
chosen in \eqref{choose} below. Since
$\partial_t\varphi+|V.\nabla_x\varphi|\leq 0$, \eqref{computation_rez} yields
$$
\Exp\Phi_N(t,\eta^N_{Nt},\hat{\eta}^N_{Nt})\leq\Exp\Phi_N(t_{n,l},\eta^N_{Nt_{n,l}},\hat{\eta}^N_{Nt_{n,l}})+t(\delta_r+\varepsilon_N+\varepsilon'_N)
$$
for any $r\in\N$ and $t\in T_{n,l}$, where
$$
\varepsilon'_N=CN^{1-d}\sum_{x\in\Omega_N}\sum_{y\not\in\Omega_N}[p(x-y)+p(y-x)]\varphi(t_{n,l},x/N)
$$
is a bound on the contribution of births and deaths in $\eta^N$. By \eqref{def_trace}, $\lim_{N\to\infty}\varepsilon'_N=0$ if $R$ and $\varepsilon$ are chosen such that
the support of $\varphi(t_{n,l},.)$ lies in $\Omega$.
%
%
Since
$$
\indicator{\{|x-x_{n,k}|\leq R-2\varepsilon-V(t-t_{n,l})\}}\leq\varphi(t,x)\leq\indicator{\{|x-x_{n,k}|\leq R-V(t-t_{n,l})\}},
$$
we obtain the  conclusion by choosing $R$ and $\varepsilon$ such that
\be\label{choose}(1+V)2^{-n}+2\varepsilon<R<2^{-n}+d(B_{n,k},\R^d\setminus\Omega)\ee
\end{proof}
\section{Proofs of Proposition \ref{prop_stat} and Theorem \ref{main_3}}
\label{proof_main_3}
%
%
%
%
%
%
%
%
%
%
A change of coordinates in \eqref{entropy_initial_vov} shows that
\begin{lemma}\label{reduce}
Let $\pi$ be a bounded Young measure on $(0,+\infty)\times(a,b)$, $n\in\R^d$, $f(.)=h(.).n$. Then $\tilde{\pi}_{t,x}:=\pi_{t,n.x}$ is a mv entropy solution
to \eqref{conservation_law} in $\Omega^n_{a,b}$ with boundary datum \eqref{boundary_datum_true_hyperplanes}, iff. $\pi$ is an entropy solution to
\be\label{one_dimension_2}\partial_t\rho(t,y)+\partial_y f(\rho(t,y))=0\ee
on $(a,b)$ with boundary datum
\be\label{boundary_datum_onedim}\lambda_a \indicator{\{a\}}+\lambda_b \indicator{\{b\}}\ee
\end{lemma}
\begin{proof}{proposition}{prop_stat}
We assume for instance $\lambda_a\leq\lambda_b$. By Lemma \ref{reduce}, we are reduced to showing that $\rho(.)$ given by \eqref{def_stat_sol} is a stationary entropy solution to \eqref{one_dimension_2} with boundary datum  \eqref{boundary_datum_onedim}.
%
%
For each $c\in\mathcal M_f(\lambda_a,\lambda_b)$, the connected component of
$\mathcal M_f(\lambda_a,\lambda_b)$ containing $c$ is a closed interval, whose left end we denote by $L(c)$.
Then $r(x):=L(\rho(x))$ is nondecreasing with respect to the usual order, thus locally of bounded variation,
and $\psi_c^\pm(r(x))$ is nonincreasing.
Hence, for $\varphi\in C^1_K(\R)$,
$$
\int_a^b \psi_c^\pm(\rho(x))\varphi'(x)dx  =  \int_a^b \psi_c^\pm(r(x))\varphi'(x)dx\geq \psi^\pm_c(r(1))\varphi(1)-\psi^\pm_c(r(0))\varphi(0)
$$
The stationary form of \eqref{entropy_initial_vov} is then a consequence of the inequality
$$
-M\phi_c^\pm(\lambda_b)\leq\psi_c^\pm(r)\leq M\phi_c^\pm(\lambda_a)
$$
which holds
for every $c\in\mathcal R$ and $r\in\mathcal M_f(\lambda_a,\lambda_b)$, provided $M\geq||f'||_\infty$, which may be assumed w.l.o.g.
\end{proof}
\mbox{}\\ \\
\subsection{Proof of Theorem \ref{main_3} for $\Omega=\Omega^n_{a,b}$}
\label{proof_main_3_1}
%
%
We first prove the following special cases of Theorem \ref{main_3}
\begin{proposition}
\label{main_3_special}
Let $\pi$  be a stationary mv entropy solution to \eqref{one_dimension_2} on $(a,b)$ with boundary datum \eqref{boundary_datum_onedim}.
Then $\pi_x=R_f(\lambda_a,\lambda_b)$ a.e. on $(a,b)$.
\end{proposition}
\begin{corollary}\label{main_3_special_bis}
Let $\rho(.,.)$  be an entropy solution to \eqref{one_dimension_2} on $(a,b)$ with boundary datum \eqref{boundary_datum_onedim}.
Then $\rho(t,.)\to R_f(\lambda_a,\lambda_b)$ in $L^1((a,b))$ as $t\to\infty$.
%
%
\end{corollary}
\begin{proof}{proposition}{main_3_special}
By \eqref{entropy_initial_vov}, $x\mapsto\pi_x(\psi_c^\pm)=:g_c^\pm(x)$
satisfies $\partial_x g^\pm_c(x)\leq 0$ in the sense of distributions on $(a,b)$,
hence it is a nonincreasing function on $(a,b)$. For any $\varphi\in C^1_K(\R)$,
$$
\int_{(a,b)}\varphi'(x)g_c^\pm(x)dx=\varphi(1)g_c^\pm(1^-)-\varphi(0)g_c^\pm(0^+)-\int_{(a,b)}\varphi(x)dg_c^\pm(x)
$$
Then, the stationary form of \eqref{entropy_initial_vov} implies
$$
g_c^\pm(0^+)\leq M\phi_c^\pm(\lambda_a),\quad g_c^\pm(1^-)\geq -M\phi_c^\pm(\lambda_b)
$$
and thus, for a.e. $x\in(a,b)$,
\be\label{equation_stat}
-M\phi_c^\pm(\lambda_b)\leq\pi_x(\psi^\pm_c)\leq M\phi^\pm_c(\lambda_a)
\ee
In \eqref{equation_stat}, we take $c=R_f(\lambda_a,\lambda_b)\in[\lambda_a,\lambda_b]$.
Since $f(\rho)-f(c)\geq 0$ for all $\rho\in\mathcal R$, the second inequality with $\phi_c^+$ and the first inequality with $\phi_c^-$ yield
$\pi_x[f(\rho)-f(c)]\leq 0$.
Hence, $\pi_x=\delta_{R_f(\lambda_a,\lambda_b)}$.
\end{proof}
\mbox{}\\ \\
For the proof of Corollary \ref{main_3_special_bis} and for subsequent use, we recall a classical contraction
and finite propagation result (see \cite{vov}).
An entropy sub-solution, resp. super-solution (\cite{bk}) to
\eqref{conservation_law}, is defined by restricting
\eqref{entropy_initial_vov} to $(\phi,\psi)=(\phi_c^+,\psi_c^+)$,
resp. $(\phi,\psi)=(\phi_c^-,\psi_c^-)$, and lies in
$C^0([0,+\infty);L^1_{\rm loc}(\Omega))$ by \cite{pan}. Mv entropy sub-solutions and super-solutions are defined in the usual way.
Remark that an entropy sub (super)-solution  remains one if the initial or boundary  datum is increased (decreased).
%
%
%
%
\begin{proposition}
\label{comparison_entropy}
Let $\pi^1(.,.)$ (resp. $\pi^2(.,.)$) be a mv entropy sub-solution
(resp. super-solution) to \eqref{conservation_law} in $\Omega$ with initial
data $\rho^1_0(.)$, $\rho^2_0(.)$ in $\Omega$ and boundary data
$\hat{\lambda}_1(.)\leq\hat{\lambda}_2(.)$ in $\partial\Omega$.
Let $C>0$ be such that $\pi^i_{t,x}$ is supported a.e. on $[0,C]$ for each $i\in\{1,2\}$.
Set $V=\sup\{|h'(\rho)|,\,\rho\leq C\}$. Then, for a.e. $t>0$ (every $t>0$ in the case of Dirac solutions),
\be \label{contraction}
\int_{\Omega\cap B(x_0,R-Vt) }\int_{{\mathcal R}^2} (\rho^1-\rho^2)^+\pi^1_{t,x}(d\rho^1)\pi^2_{t,x}(d\rho^2) dx
\leq
\int_{\Omega\cap B(x_0,R)}(\rho^1_0(x)-\rho^2_0(x))^+ dx \ee
\end{proposition}
In particular, \eqref{contraction} implies a maximum principle for $\pi^1$, (and a similar minimum principle for $\pi^2$),
$\pi^1_{t,x}([0,C_1])=1$, where
$$
C_1:=
\max\left(
\sup_\Omega(\rho^1_0(.)),\sup_{\partial\Omega}\hat{\lambda}_1(.)
\right)
$$
%
%
%
\begin{proof}{corollary}{main_3_special_bis}
Let $\delta_1$ be a smooth nonnegative function, supported on $[0,1]$, such that $\int_0^1\delta_1(u)=1$.
For $s>0$ and $x\in(a,b)$, let $\chi^{s,+}(t)=s^{-1}\delta_1[(t-s)/s]$, $\chi^{s,-}(t)=2s^{-1}\delta_1[(2t-s)/s]$, and

$$\pi^{s,\pm}_x=\int_0^{+\infty}\delta_{\rho(t,x)}\chi^{s,\pm}(t)dt$$
The Young measures $\pi^{s,\pm}$ are uniformly bounded.
Taking $\varphi(t,x)=\varphi_0(x)\chi^{s,\pm}(t)$ in \eqref{entropy_initial_vov}, with arbitrary $\varphi_0\in C^1_K(\R)$,
shows that any subsequential limit $\pi$ of $\pi^{s,\pm}$ as $s\to\infty$ is a stationary
mv entropy solution to \eqref{one_dimension_2} on $(a,b)$ with boundary datum \eqref{boundary_datum_onedim}. Hence, by Proposition \ref{main_3_special}, $\pi=\delta_{R_f(\lambda_a,\lambda_b)}$. Thus,
\begin{eqnarray*}
\lim_{s\to\infty}\int_{0}^{+\infty}\int_a^b|\rho(t,x)-R_f(\lambda_a,\lambda_b)|\chi^{s,\pm}(t)dxdt & = & 0 \\
\end{eqnarray*}
Let ${\mathcal I}^{s,\pm}$ denote the above integral. Since the constant $R_f(\lambda_a,\lambda_b)$ is a stationary entropy solution to \eqref{one_dimension_2} on $(a,b)$ with boundary datum \eqref{boundary_datum_onedim}, and $\chi^{s,+}$, $\chi^{s,-}$ are supported respectively on $[s,2s]$ and $[s/2,s]$, by Lemma \ref{lemma_comparison},
$$
\mathcal I^{s,+}\leq\int_a^b|\rho(s,x)-R_f(\lambda_a,\lambda_b)|dx\leq\mathcal I^{s,-}
$$
whence the result follows.
\end{proof}
\mbox{}\\ \\
\begin{proof}{theorem}{main_3}
{\em Part one.}
Let $C>0$ be such that $\pi_{t,x}$ is supported on $[0,C]$ and $\rho(t,x)\leq C$ for a.e. $(t,x)$. For $r\in\{0,C\}$, we denote by $\tilde{\rho}^r(t,x)$ the entropy solution
to \eqref{conservation_law} in $\Omega^n_{a,b}$ with uniform initial datum $r$, and boundary datum \eqref{boundary_datum_true_hyperplanes}.
Since these data are invariant by translations orthogonal to $n$, the same holds for $\tilde{\rho}^r$. Hence, by Lemma \ref{reduce}, $\tilde{\rho}^r(t,x)=\rho^r(t,n.x)$, where $\rho^r$ is the entropy solution to \eqref{one_dimension_2} in $(a,b)$, with uniform initial datum $r$ and boundary datum \eqref{boundary_datum_onedim}.
On the other hand, $\pi_x$ and $\delta_{\rho(t,x)}$ are mv entropy sub-solutions (resp. super-solutions) for the boundary datum \eqref{boundary_datum_true_hyperplanes} and
uniform initial datum $C$ (resp. $0$).  Thus, by Proposition \ref{lemma_comparison}, $\tilde{\rho}^0(t,x)\leq\rho(t,x)\leq\tilde{\rho}^C(t,x)$, and
$\pi_x([\tilde{\rho}^0(t,x),\tilde{\rho}^C(t,x)])=1$, for a.e. $(t,x)$. Conclusions follow from Corollary \ref{main_3_special_bis} applied to ${\rho}^r$.
\end{proof}
\subsection{Proof of Theorem \ref{main_3} for $\Omega$ satisfying \eqref{perturbation}}
\label{proof_main_3_2}
We shall need the following
\begin{lemma}
\label{lemma_comparison} Let $\rho(.,.)$ be the entropy solution to
\eqref{conservation_law} in $\Omega_1\subset\R^d$ with data
$\rho_0(.)$ in $\Omega_1$, $\hat{\lambda}(.)$ on $\partial\Omega_1$.
Let $\Omega_2$ be an open subset of $\Omega_1$ with locally finite perimeter,
$\Sigma\subset\partial\Omega_1\cap\partial\Omega_2$, and $R^\pm$
constants such that $R^-\leq\rho\leq R^+$ a.e. in
$(0,+\infty)\times\Omega_2$. Let $\hat{\lambda}^\pm(.)$ be equal to
$\hat{\lambda}(.)$ on $\Sigma$ and $R^\pm$ on
$\partial\Omega_2\backslash\Sigma$.
%
%
Then $\rho(.,.)$ is an entropy sub-solution (resp. super-solution)
to \eqref{conservation_law} in $\Omega_2$ for data $\rho_0(.)$
restricted to $\Omega_2$ and $\hat{\lambda}^+(.)$ (resp.
$\hat{\lambda}^-$(.)) on $\partial\Omega_2$ .
\end{lemma}
\begin{proof}{lemma}{lemma_comparison}
We may assume w.l.o.g. that $M$ in \eqref{entropy_initial_vov} is a Lispchitz constant for $h$ on $[0,||\rho||_\infty]$. Let $(\phi^\pm,\psi^\pm)=(\phi^\pm_c,\psi^\pm_c)$. By
\eqref{entropy_initial_vov},
$m^\pm(dt,dx):=\partial_t\phi^\pm(\rho)+\nabla_x.\psi^\pm(\rho)
$
is a nonpositive measure
on $(0,+\infty)\times\Omega_1$, hence $(\phi^\pm(\rho),\psi^\pm(\rho))$ is a divergence-measure field. Let ${\mathcal O}$ be an open subset
of $\Omega_1$ with locally finite perimeter. By the generalized Gauss-Green formula
(\cite{cf, ct}),
\begin{eqnarray} \nonumber
I_{\mathcal O}(\varphi) & := & \int_{(0,+\infty)\times{\mathcal
O}}[\phi^\pm(\rho(t,x))\partial_t\varphi(t,x)  + \psi^\pm(\rho(t,x)).\nabla_x\varphi(t,x)
]dt dx \\
& = &  -\int_{(0,+\infty)\times\partial{\mathcal
O}}\varphi(t,x)\hat{\psi}^\pm(t,x)dt \,d{\mathcal
H}^{d-1}(x)-\int_{\mathcal O}\varphi(0,x){\phi}^\pm(\rho_0(x))dx\nonumber\\
& & -\int_{(0,+\infty)\times{\mathcal
O}}\varphi(t,x)m^\pm(dt,dx)\label{gauss_green_1}
\end{eqnarray}
for every $0\leq \varphi\in C^1_K([0,+\infty)\times\R^d)$, where
$\hat{\psi}^\pm(t,x)$
is a weak normal trace for $\psi^\pm(\rho(t,x))$ on $\partial\mathcal O$.
For ${\mathcal O}=\Omega_1$, \eqref{entropy_initial_vov} implies
\be \label{trace_omega}
M\phi^\pm(\hat{\lambda}(x))-\hat{\psi}^\pm(t,x)\geq 0 \ee
a.e. on $\partial\Omega_1$. To evaluate the l.h.s. of \eqref{entropy_initial_vov} for $\Omega=\Omega_2$, we let ${\mathcal O}=\Omega_2$ in \eqref{gauss_green_1}.
We obtain boundary integrands \eqref{trace_omega} on $\Sigma$, and
$M\phi^\pm(R^\pm)-\hat{\psi}^\pm(t,x)$ on $\partial\Omega_2\setminus\Sigma$.
%
%
The latter is nonnegative because,
for a.e. $(t,x)\in(0,+\infty)\times\Omega_2$ and any unitary vector
$n\in\R^d$,
$\psi^{\pm}(\rho(t,x)).n\leq
M\phi^\pm(\rho(t,x))\leq M\phi^\pm(R^\pm)$.
\end{proof}
\mbox{}\\ \\
\begin{lemma}
\label{flat_general} Let $\rho(.,.)$ be the entropy solution to
\eqref{conservation_law} with uniform data $r\in\mathcal R$ in
$\Omega$ and $\lambda\in{\mathcal R}$ on $\partial\Omega$. Then $\rho(t,.)\to\lambda$ as $t\to\infty$ in
$L^1_{\rm loc}(\overline{\Omega})$.
\end{lemma}
\begin{proof}{lemma}{flat_general}
Denote by $\rho_{a',b'}(t,x)$ the entropy solution to
\eqref{conservation_law} in $\Omega^n_{a',b'}$ with uniform data $r$
in $\Omega^n_{a',b'}$ and $\lambda$ on $\partial\Omega^n_{a',b'}$.
By Lemma \ref{reduce} and Corollary \ref{main_3_special_bis},
\be \label{bylemmas} \rho_{a',b'}(t,.)\to\lambda,\quad\mbox{in }
L^1_{\rm loc}(\overline{\Omega}^n_{a',b'}) \ee
We consider the case $r\geq\lambda$, the case $r\leq\lambda$ being
similar. By maximum principle,
\be\label{yield_1}\rho(t,.)\geq\lambda, \quad
\rho_{a',b'}(t,.)\geq\lambda\ee
for every $t>0$. Hence, by Lemma \ref{lemma_comparison}, the
restriction of $\rho_{a',b'}$ to $(0,+\infty)\times\Omega$ is an
entropy super-solution to \eqref{conservation_law} in $\Omega$ for
the uniform data $r$ in $\Omega$, $\lambda$ on $\partial\Omega$.
Thus, by Proposition \ref{comparison_entropy},
$\rho_{a',b'}(t,.)\geq\rho(t,.)$ in $\Omega$ for every $t>0$. This,
\eqref{yield_1} and \eqref{bylemmas} imply the result.
\end{proof}
\mbox{}\\ \\
\begin{proof}{Theorem}{main_3}
{\em Part two.}
We first prove (ii).
We consider $\lambda_a\leq\lambda_b$, the
reverse case being similar. We set $f(.)=h(.).n$ and denote by $T_t(\Omega,\hat{\lambda})$
the solution semigroup for \eqref{conservation_law} in $\Omega$ with
boundary datum $\hat{\lambda}(.)$, i.e. $T_t(\Omega,\hat{\lambda})\rho_0(x)$ is the entropy solution at time $t$ with initial datum $\rho_0(.)$.
Proposition \ref{comparison_entropy} implies that $T_t$ is monotone.
We define the following entropy solutions to
\eqref{conservation_law}. For $\gamma\in\{a,b\}$,
\begin{eqnarray*} \rho^\gamma(t,.)=T_t(\Omega,\hat{\lambda})[\lambda^f_\gamma\indicator{\Omega}] & , &
\overline{\rho}^\gamma(t,.)=T_t(\Omega,\lambda^f_\gamma\indicator{\partial\Omega})\rho_0(.)\\
{\rho}^{a',b}(t,.)=T_t(\Omega^n_{a',b},\hat{\lambda}^{a',b})[\lambda^f_b\indicator{\Omega^n_{a',b}}] & , &
{\rho}^{a,b'}(t,.)=T_t(\Omega^n_{a,b'},\hat{\lambda}^{a,b'})[\lambda^f_a\indicator{\Omega^n_{a,b'}}]
\end{eqnarray*}
where
\begin{eqnarray*} \hat{\lambda}^{a',b}(x)  := \lambda_a\indicator{\{n.x=a'\}}+\lambda^f_b\indicator{\{n.x=b\}} & , &
%
\hat{\lambda}^{a,b'}(x) :=  \lambda^f_a\indicator{\{n.x=a\}}+\lambda_b\indicator{\{n.x=b'\}}
\end{eqnarray*}
%
%
By maximum principle,
\begin{eqnarray*} \lambda^{f}_a & \leq &
\min(\rho^{a}(t,.),\,{\rho}^{a,b'}(t,.))\leq
\max(\rho^a(t,.),\,{\rho}^{a,b'}(t,.))\leq\lambda_b\\
\lambda_a & \leq &
\min(\rho^b(t,.),\,{\rho}^{a',b}(t,.))\leq
\max(\rho^b(t,.),\,{\rho}^{a',b}(t,.))\leq\lambda^f_b
\end{eqnarray*}
for every $t>0$, respectively  on $\Omega\cap\Omega^n_{a,b'}$ and $\Omega\cap\Omega^n_{a',b}$.
It follows by Lemma \ref{lemma_comparison} that restrictions of ${\rho}^a$ and ${\rho}^{a,b'}$
to $\Omega\cap\Omega^n_{a,b'}$ are respectively a super-solution and
a sub-solution to \eqref{conservation_law} in
$\Omega\cap\Omega^n_{a,b'}$ for the uniform initial datum
$\lambda^{f}_a$, and boundary datum
$$\hat{\lambda}^a(x)=\lambda^{f}_a\indicator{\{n.x=a\}}+\lambda_b\indicator{\partial\Omega_b}$$
Similarly, restrictions of ${\rho}^b$ and ${\rho}^{a',b}$ to
$\Omega\cap\Omega^n_{a',b}$ are respectively a sub-solution and a
super-solution to \eqref{conservation_law} in
$\Omega\cap\Omega^n_{a',b}$ for the uniform initial datum
$\lambda^{f}_b$, and boundary datum
$$\hat{\lambda}^b(x)=
%
%
\lambda_a\indicator{\partial\Omega_a}+\lambda^f_b\indicator
{\{n.x=b\}}$$
%
%
%
%
%
Let $0<s<t$ and $B(x_0,\delta)\subset\Omega$.
By repeated use of Proposition \ref{comparison_entropy}, we obtain first
$\lambda^f_a\leq\overline{\rho}^a\leq\rho\leq\overline{\rho}^b\leq\lambda^f_b$ on $\Omega$, and then
\begin{eqnarray}
& \dsp\int_{\Omega\cap\Omega^n_{a',b}\cap
B(x_0,\delta)}(\rho(t,x)-{\rho}^{a',b}(s,x))^+dx
& \label{epsilon_5}\\
%
%
\leq & \dsp\int_{\Omega\cap
B(x_0,\delta)}(T_s(\Omega,\hat{\lambda})\overline{\rho}^b(t-s,.)(x)-T_s(\Omega,\hat{\lambda})[\lambda^f_b\indicator{\Omega}](x))^+
dx & \nonumber\\
+ &
\dsp\int_{\Omega\cap \Omega^n_{a',b}\cap
B(x_0,\delta)}(\rho^b(s,x)-{\rho}^{a',b}(s,x))^+
dx \leq \int_{\Omega\cap
B(x_0,\delta+Vs)}(\overline{\rho}^b(t-s,x)-\lambda^b_f)^+ &
dx\nonumber
%
%
\end{eqnarray}
By Lemma \ref{reduce}, Corollary \ref{main_3_special_bis} and Lemma \ref{flat_general},
${\rho}^{a',b}(t,.)\to R_{f}(\lambda_a,\lambda_b)$ in
$L^1_{\rm loc}(\overline{\Omega}^n_{a',b})$, and  $\overline{\rho}^b(t,.)\to \lambda^f_b$ in
$L^1_{\rm loc}(\overline{\Omega})$ as $t\to\infty$.
Thus, $t\to\infty$ and $s\to\infty$ in \eqref{epsilon_5} yields
$$\lim_{t\to\infty}\int_{\Omega\cap\Omega^n_{a',b}\cap B(x_0,\delta)}(\rho(t,x)-R_f(\lambda_a,\lambda_b))^+dx=0$$
The  negative part on $\Omega\cap\Omega^n_{a,b'}$ is treated in a similar way using
$\rho^{a,b'}$ and $\overline{\rho}^a$. Note that $R_f(\lambda_a,\lambda_b)=\lambda_\gamma$, with $\gamma\in\{a,b\}$, necessarily implies $\lambda_\gamma=\lambda^f_\gamma$.
%
%
Statement (i) follows from (ii)  as in the first part of the proof
in Subsection \ref{proof_main_3_1}, by comparing $\pi$ with entropy solutions with uniform initial datum in $\Omega$.
\end{proof}
\begin{appendix}
\section{Exclusion process with overtaking}
\label{appendix_kstep}
The following model is closely related to the $k$-step exclusion process introduced in \cite{gui}.
\subsection{The model}\label{model_kstep}
Let $\mathcal K=\mathcal R=1$, $\left(e_1,\ldots,e_d\right)$ be the canonical
basis of $\R^d$, ${\mathcal D}=\{\pm e_1,\ldots,\pm e_d\}$,
$k\in\N$, and $(\beta^\alpha_j)_{\alpha\in\mathcal D,\,j\in\N}$ be a family of
nonnegative real numbers such that
\be\label{summable_beta}
\sum_{j\in\N}j\beta_j^\alpha<+\infty,\quad\forall\alpha\in\mathcal D
\ee
%
%
%
Thanks to \eqref{summable_beta}, following the lines of \cite{lig}, one can define a Feller process on $\cspace$ with  generator
\be \label{gen_kstep} \gen f(\eta)=\sum_{x\in\Z^d}\sum_{\alpha\in{\mathcal D}}\sum_{j\in\N}
\beta_j^\alpha c^\alpha_{x,j}(\eta)\left[f\left(\eta^{x,x+j\alpha}\right)-f(\eta)\right] \ee
where \be \label{def_ksteprates} c^\alpha_{x,j}(\eta):=
(1-\eta(x+j\alpha))\prod_{n=0}^{j-1}\eta(x+n\alpha) \ee
The interpretation is that  a particle jumps to the first vacant site in a randomly chosen
direction $\alpha=\pm e_i$ with rate $\beta^\alpha_j$, where $j$ is the distance to the target site, or stays where it is if no vacant site is found.
For every $\rho\in[0,1]$, the product Bernoulli measure $\nu_\rho$, whose marginal $\theta_\rho$ at each site is the Bernoulli measure with parameter $\rho$, is invariant for this
process.
We further make the irreducibility assumption
$\beta^\alpha_1+\beta^{-\alpha}_{1}>0$ and monotonicity assumption $\beta^\alpha_{j+1}\leq\beta^\alpha_j$ for all $\alpha\in\mathcal D$, the latter implying attractiveness.
These two assumptions imply  \eqref{extremal_measures} and Proposition \ref{invariant_coupled}, which can be established
as in \cite{gui}. \\ \\
\textbf{Open-boundary dynamics.}
For this model, \eqref{def_genopen} can be expressed explicitely as follows.
Given $x\in \Z^d$, $\alpha\in\mathcal D$ and $j\in\N$, let
\be
\overline{c}^\alpha_{x,j}(\eta):=\prod_{z\in[x,x+j\alpha]\cap\Omega_N}\eta(z)
\prod_{z\in[x,x+j\alpha]\backslash\Omega_N}\lambda_N(z)
\ee
where $[x,x+j\alpha]:=\{x+i\alpha:\,i=0,\ldots,j\}$. If $x\in\Omega_N$ and $x+j\alpha\in\Omega_N$, a jump from $x$ to $x+j\alpha$ occurs at rate
$\overline{c}^\alpha_{x,j}$. Note that this rate may depend on the reservoir profile if $[x,x+j\alpha]\not\subset\Omega_N$.
If $x\in\Z^d\setminus\Omega_N$ and $x+j\alpha\in\Omega_N$, a particle is created at $x+j\alpha$ with rate $\overline{c}^\alpha_{x,j}$.
If $x+j\alpha\in\Z^d\backslash\Omega_N$ and $x\in\Omega_N$, a particle is removed from $x$ at rate $\overline{c}^\alpha_{x,j}$.\\ \\
%
%
\emph{Example.} Let $d=1$, $\Omega=(0,1)$, $\Omega_N=\{1,\ldots,N-1\}$,  $\mathcal D=\{-1,1\}$, $\beta_1^{-1}=0$ (no jumps to the left),
$\beta_1^1=\beta_1\geq\beta_2^1=\beta_2>0=\beta^1_3$, $\lambda_N(-1)=\lambda_N(0)=\lambda_a$, $\lambda_N(N)=\lambda_N(N+1)=\lambda_a$,
$\hat{\lambda}=\lambda_a\indicator{\{0\}}+\lambda_b\indicator{\{1\}}$.
The boundary dynamics is defined as follows. A particle is created
at site $1$ with rate
$(\beta_1\lambda_l+\beta_2\lambda_l^2)[1-\eta(1)]$ and at site $2$
with rate $\beta_2\lambda_l\eta(1)[1-\eta(2)]$, removed from site
$N-2$ with rate $\beta_2\eta(N-2)\eta(N-1)[1-\lambda_r]$ and from
site $N-1$ with rate
$[\beta_1(1-\lambda_r)+\beta_2\lambda_r(1-\lambda_r)]\eta(N-1)$.\\ \\
%
%
\textbf{Flux function.}
The microscopic flux function involved in \eqref{current_macro} is now
\be j(\eta)  =  \sum_{\alpha\in\mathcal
D}\sum_{j\in\N}j\alpha\beta^\alpha_j c^\alpha_{0,j}(\eta),
\label{current_kstep} \ee
which results in
\be\label{flux_kstep}
h(\rho)=\rho(1-\rho)\sum_{i=1}^d \left[
\sum_{j\in\N}j\left(\beta_j^{e_i}-\beta_j^{-e_i}\right)\rho^{j-1}
\right]e_i
\ee
Assumption \eqref{summable_beta}
%
%
%
implies $h\in C^1([0,1])$. In particular, any flux function of the form $h(\rho)=\rho(1-\rho)P(\rho)$, where $P$ is a $\R^d$-valued polynomial, can be obtained by a suitable choice of $\beta^\alpha_j$.
\subsection{Proof of Theorems \ref{main} and \ref{main_2}}
For this model, the range $\Delta$ of the dynamics involved in \eqref{def_trace} is now defined by $\Delta=\sup\{j\in\N:\,\sum_{\alpha\in\mathcal D}\beta^\alpha_j>0\}$.
The coupled generator for the process on $\Z^d$ is (see \cite{gui} for a similar coupling), with $(\eta,\xi)\in\cspace^2$,
\begin{eqnarray*}
\widetilde{L}f(\eta,\xi) & = &  \sum_{\alpha\in\mathcal D}\sum_{x\in\Z^d}\sum_{k,l\in\N}
\beta_k^\alpha\wedge\beta_l^\alpha {c}^\alpha_{x,k}(\eta){c}^\alpha_{x,l}(\xi)\left[f\left(\eta^{x,x+k},\xi^{x,x+l}\right)-f(\eta,\xi)\right]\nonumber\\
& + & \sum_{\alpha\in\mathcal D}\sum_{x\in\Z^d}\sum_{k,l\in\N}
(\beta_k^\alpha-\beta_l^\alpha)^+ {c}^\alpha_{x,k}(\eta){c}^\alpha_{x,l}(\xi)\left[f\left(\eta^{x,x+k},\xi\right)-f(\eta,\xi)\right]\nonumber \\
& + & \sum_{\alpha\in\mathcal D}\sum_{x\in\Z^d}\sum_{k,l\in\N}
(\beta_l^\alpha-\beta_k^\alpha)^+ {c}^\alpha_{x,k}(\eta){c}^\alpha_{x,l}(\xi)\left[f\left(\eta,\xi^{x,x+l}\right)-f(\eta,\xi)\right]\\
& + & \sum_{\alpha\in\mathcal D}\sum_{x\in\Z^d}\sum_{k\in\N}
\beta_k^\alpha {c}^\alpha_{x,k}(\eta)\left[1-\sum_{l\in\N}{c}^\alpha_{x,l}(\xi)\right]\left[f\left(\eta^{x,x+k},\xi\right)-f(\eta,\xi)\right]\nonumber \\
& + & \sum_{\alpha\in\mathcal D}\sum_{x\in\Z^d}\sum_{l\in\N}
\beta_l^\alpha\left[1-
\sum_{k\in\N}{c}^\alpha_{x,k}(\eta)\right]{c}^\alpha_{x,l}(\xi)\left[f\left(\eta,\xi^{x,x+l}\right)-f(\eta,\xi)\right]
\end{eqnarray*}
%
%
This coupling preserves order as a result of the monotonicity assumption on $\beta_j^\alpha$.
For the open coupled process \eqref{def_couplegen}, the reservoir measure \eqref{outer_product_tilde} is defined here from Bernoulli marginals \eqref{bernoulli}.
%
%
Lemma \ref{equilibrium_property}, which depends only on product invariant measures and model-independent definition \eqref{def_genopen}, still holds true.
For Lemma \ref{monotone_coupling}, a direct proof, based on explicit coupling rates for the open process, was given for the Misanthrope's process. Such a proof is tedious here.
However, Lemma \ref{monotone_coupling} holds more generally because the above coupling on $\Z^d$ is monotone, and the model-independent definition \eqref{def_couplegen} automatically inherits this property:\\
\begin{proof}{lemma}{monotone_coupling}
{\em Model-independent version.}
We have to show that whenever $\eta\leq\xi$, where $(\eta,\xi)\in\cspace_N^2$, a transition $(\eta,\xi)\to(\eta',\xi')$ in the coupled process \eqref{def_couplegen} cannot have positive rate unless $\eta'\leq\xi'$.
From \eqref{def_couplegen}, the rate of this transition is given by
\be\label{coupled_transition}
\widetilde{C}_N(\eta,\xi;\eta',\xi')=\int_{\overline{\cspace}_N^2}\sum_{(\overline{\eta'},\overline{\xi'})\in\overline{\cspace}_N^2}\widetilde{C}(\eta\oplus\overline{\eta},\xi\oplus\overline{\xi};\eta'\oplus\overline{\eta'},\xi'\oplus\overline{\xi'})
d\widetilde{\overline{\nu}}_N(\overline{\eta},\overline{\xi})\ee
where the integrand $\widetilde{C}$ denotes the rate of a transition
$$
(\eta\oplus\overline{\eta},\xi\oplus\overline{\xi})\to(\eta'\oplus\overline{\eta'},\xi'\oplus\overline{\xi'})
$$
for the  coupled generator $\widetilde{\gen}$ on $\Z^d$.
Since $\lambda_N(x)\leq c$ for every $x\in\Z^d\backslash\Omega_N$, $\widetilde{\overline{\nu}}_N$ is supported on pairs $(\overline{\eta},\overline{\xi})$ such that $\overline{\eta}\leq\overline{\xi}$. For such pairs, the rate \eqref{coupled_transition} is zero if $\eta'\not\leq\xi'$, because
$\widetilde{\gen}$ is a  monotone coupling.
\end{proof}
\mbox{}\\ \\
We are left to establish analogues of \eqref{replacement_boundary}--\eqref{replacement_bulk} for this model.
For simplicity we will only consider $\Phi^+_N$.
For the process on $\Z^d$ we have, with $(\eta,\xi)\in\cspace^2$,
\begin{eqnarray}
\widetilde{\gen}\Phi^+_N(t,\eta,\xi) & = & \sum_{\alpha\in\mathcal
D}\sum_{x\in\Z^d}\sum_{k,l\in\N:\,k> l}\beta^\alpha_k
c_{x,k}^\alpha(\eta)
c^\alpha_{x,l}(\xi)[1-\xi(x+k\alpha)]\nabla^\alpha_{N,l,k}\left[\varphi\indicator{\Omega_N}\right](t,x)
\nonumber\\
& + & \sum_{\alpha\in\mathcal D}\sum_{x\in\Z^d}\sum_{k,l\in\N:\,k>
l}[\beta^\alpha_l-\beta^\alpha_k]c_{x,k}^\alpha(\eta)c^\alpha_{x,l}(\xi)
\nabla^\alpha_{N,l,0}\left[\varphi\indicator{\Omega_N}\right](t,x)
%
\nonumber\\
& + & \sum_{\alpha\in\mathcal
D}\sum_{x\in\Z^d}\sum_{k\in\N}\beta^\alpha_k
c_{x,k}^\alpha(\eta)[1-\xi(x)][1-\xi(x+k\alpha)]
\nabla^\alpha_{N,0,k}\left[\varphi\indicator{\Omega_N}\right](t,x)
\nonumber\\
& - & \sum_{\alpha\in\mathcal D}\sum_{x\in\Z^d}\sum_{k,l\in\N:\,k\neq l}\beta^\alpha_k\wedge\beta^\alpha_l c_{x,k}^\alpha(\eta)c^\alpha_{x,l}(\xi)\xi(x+k\alpha)
\left[\varphi\indicator{\Omega_N}\right]\left(t,\frac{x+l\alpha}{N}\right)\nonumber\\
& - & \sum_{\alpha\in\mathcal D}\sum_{x\in\Z^d}\sum_{k\in\N}\beta^\alpha_k c_{x,k}^\alpha(\eta)[1-\xi(x)]\xi(x+k\alpha)\left[\varphi\indicator{\Omega_N}\right]\left(t,\frac{x}{N}\right)\nonumber\\
& - & \sum_{\alpha\in\mathcal
D}\sum_{x\in\Z^d}\sum_{l\in\N}\beta^\alpha_l
c_{x,l}^\alpha(\xi)[1-\eta(x)]\eta(x+l\alpha)\left[\varphi\indicator{\Omega_N}\right]\left(t,\frac{x+l\alpha}{N}\right)\label{genphi_kstep}
\end{eqnarray}
where $\nabla^\alpha_{N,l,k}\varphi(t,x):=\varphi[t,(x+k\alpha)/N]-\varphi[t,(x+l\alpha)/N]$. An upper bound for  \eqref{genphi_kstep} is obtained by gathering
three types of terms on the first three lines:\\ \\
\textbf{Boundary terms} are those whose last factor (the variation of $\Phi_N^+$) is positive of order $1$:
(i) those on the first line for which $x+k\alpha\in\Omega_N,\,x+l\alpha\not\in\Omega_N$;
(ii) those on the second line for which $x\in\Omega_N,\,x+l\alpha\not\in\Omega_N$;
(iii) those on the third line for which $x\not\in\Omega_N,\,x+k\alpha\in\Omega_N$.
The corresponding variation of $\Phi_N^+$ is $\varphi(z/N)$, where $z$ is the site inside $\Omega_N$ in (i)--(iii).
We denote by $\widetilde{B}^i_N(t,\eta,\xi)$ the sum of boundary terms from the $i$-th line of \eqref{genphi_kstep}.
For $(\eta,\xi)\in\cspace_N^2$, we denote by
$$
\widetilde{B}^i_{N,c}(t,\eta,\xi)=\int_{\overline{\cspace}_N^2}\widetilde{B}^i_N(t,\eta\oplus\overline{\eta},\xi\oplus\overline{\xi})d\widetilde{\overline{\nu}}_N(\overline{\eta},\overline{\xi})
$$
the resulting contribution to $\widetilde{L}_{N,c}\Phi^+_N$. This averaging produces a factor $(\lambda_N(x+l)-c)^+$ on the first two lines,
and a factor $(\lambda_N(x)-c)^+$ on the third line of \eqref{genphi_kstep}. After changes of index and exchanges of summations, we obtain the bounds
\begin{eqnarray*}
\widetilde{B}^{i}_{N,c}(t,\eta,\xi) & \leq & \sum_{\alpha\in\mathcal D}\sum_{k=1}^{+\infty}\beta_k^\alpha\sum_{u\not\in\Omega_N,\,d(u,\Omega_N)\leq k}(\lambda_N(u)-c)^+\sup_{v\in\Omega_N,\,|v-u|\leq k}\varphi\left(\frac{v}{N}\right)
%
%
\end{eqnarray*}
By \eqref{def_trace} and \eqref{summable_beta}, the above r.h.s. is bounded by a constant times $\mathcal K^\pm_{\partial\Omega}(\varphi)$.
Thus the contribution of boundary terms produces the r.h.s. of \eqref{replacement_boundary}.\\ \\
\textbf{Bulk terms} are those
%
for which $[x,x+k\alpha]\subset\Omega_N$.
We introduce a cutoff $k\leq r$, similar to \eqref{current_plus}--\eqref{current_min}. By \eqref{def_trace} and \eqref{summable_beta}, the contribution of terms $k>r$ yields a vanishing $\delta_r$ as in \eqref{computation_rez}. The contribution of terms $k\leq r$ is analogous to the leading term on the r.h.s. of \eqref{computation_rez},
with \eqref{current_plus} replaced by
\begin{eqnarray*}
\tilde{j}_r^{+}(\eta,\xi) & = & \sum_{\alpha\in\mathcal D}\sum_{k,l\in\N:\,r\geq k>l}(k-l)\beta^\alpha_k c^\alpha_{0,k}(\eta)c^\alpha_{0,l}(\xi)[1-\xi(k\alpha)]\\
& - & \sum_{\alpha\in\mathcal D}\sum_{k,l\in\N:\,r\geq k>l}l[\beta^\alpha_l-\beta^\alpha_k] c^\alpha_{0,k}(\eta)c^\alpha_{0,l}(\xi)\\
& + & \sum_{\alpha\in\mathcal D}\sum_{k\in\N:\,r\geq k}k\beta^\alpha_k c^\alpha_{0,k}(\eta)[1-\xi(0)][1-\xi(k\alpha)]
%
%
%
%
%
%
\end{eqnarray*}
Let $\tilde{j}^+$ be defined as $\tilde{j}^+_r$ without the truncation $k\leq r$.
%
%
Note that $\tilde{j}^+(\eta,\xi)=0$ if $\eta\leq\xi$, and $\tilde{j}^+(\eta,\xi)=j(\eta)-j(\xi)$ if $\eta\geq\xi$, with $j$
given by \eqref{current_kstep}.
Hence, by Proposition \ref{invariant_coupled}, $<\tilde{j}^+>_{\rho,c}=\psi^+_c(\rho)$. By \eqref{summable_beta},
$\tilde{j}_r^+\to\tilde{j}^+$ uniformly as $r\to\infty$.
Thus, by Proposition \ref{loc_eq}, the contribution of bulk terms
produces the r.h.s. of \eqref{replacement_bulk}. \\ \\
\textbf{Error terms} are those not considered yet,
%
for which $[x,x+k\alpha]\not\subset\Omega_N$ and $[x,x+k\alpha]\cap\Omega_N\neq\emptyset$.
For these terms, $d(x,\Omega_N)\leq k$ and $d(x,\Z^d\setminus\Omega_N)\leq k$, and the corresponding gradient term in \eqref{genphi_kstep}
is or order $N^{-1}$. By \eqref{def_trace} and \eqref{summable_beta}, the total contribution of such terms is $O(N^{-1})$.
\end{appendix}
\mbox{}\\ \\
\textbf{Acknowledgements.} I thank Organizers of the 2004 Oberwolfach Workshop ``Large scale stochastic dynamics'' (C. Landim, S. Olla, H. Spohn), as well as T. Bodineau, R.
Esposito, S. Grosskinsky, J.L. Lebowitz, R. Marra, G. Sch\"{u}tz and E.R. Speer for invitations and discussions.
I thank anonymous referees for useful suggestions to improve the presentation of the paper.
%
%

%
\end{document}